\documentclass[a4paper,twoside, 11pt]{article}
\usepackage{amssymb,amsfonts,amsmath}
\usepackage{fancyhdr} 

\newtheorem{defin}{~~~~Definition}
\newtheorem{prop}{~~~~Proposition}[section]
\newtheorem{remark}{~~~~Remark}[section]
\newtheorem{theor}{~~~~Theorem}
\newtheorem{lemma}{~~~~Lemma}[section]

\newcommand{\vf}{\varphi}
\newcommand{\om}{\omega}
\newcommand{\bom}{\overline \omega}

\title{Fundamental form and Cartan's tensor of (2,5)-distributions coincide} 
\date{}
\author{Igor Zelenko\thanks{S.I.S.S.A., Via Beirut 2-4, 34014, Trieste, Italy; email: zelenko@sissa.it}} 

\begin{document}
	
	\maketitle
	\setcounter{equation}{0} 
	
	\begin{abstract}
		In our previous paper \cite{zelduke} for generic rank 2 vector
		distributions on $n$-dimensional manifold ($n\geq 5$) we   
		constructed  a special differential invariant, the {\it  
			fundamental form}. In the case $n=5$ this differential invariant
		has the same algebraic nature, 
		as the {\it covariant binary biquadratic form}, constructed 
		by E.Cartan in \cite{cartan}, using his ``reduction- 
		prolongation'' procedure (we call this form 
		Cartan's tensor). 
		In the present paper 
		we prove that our  fundamental form coincides   
		(up to a constant factor $-35$) with Cartan's tensor. This result explains
		geometric reason for  existence of Cartan's tensor (originally this tensor was 
		obtained by very sophisticated algebraic manipulations) and 
		gives the true analogs of this tensor in Riemannian 
		geometry. In addition, as a part of the proof, we obtain a 
		new useful formula for Cartan's tensor in terms of 
		structural functions of any frame naturally adapted to the 
		distribution. 
	\end{abstract}
	

\section{Introduction}
\setcounter{equation}{0}
 \indent 

 A rank $l$ vector distribution $D$ on the $n$-dimensional
manifold $M$ or $(l,n)$-distribution (where $l<n$) is by 
definition an $l$-dimensional subbundle of the tangent 
bundle $TM$. 
In other words, for each point $q\in M$ a $l$-dimensional 
subspace $D(q)$ of the tangent space $T_qM$ is chosen and 
$D(q)$ depends smoothly on $q$. Two vector distributions 
$D_1$ and $D_2$ are called equivalent, if there exists a 
diffeomorphism $F:M\mapsto M$ such that 
$F_*D_1(q)=D_2(F(q))$ for any $q\in M$. Two germs of vector 
distributions $D_1$ and $D_2$ at the point $q_0\in M$ are 
called equivalent, if there exist neighborhoods $U$ and 
$\tilde U$ of $q_0$ and a diffeomorphism $F:U\mapsto \tilde 
U$ such that 
$$\begin{array}{c}F_*D_1(q)=D_2(F(q)),\,\,\forall q\in 
U;\\F(q_0)=q_0.\end{array}$$

An obvious (but very rough in the most cases) invariant of 
distribution $D$ at $q$ is so-called {\it small growth 
vector} at $q$: it is the tuple $$\bigl(\dim D(q),\dim 
D^2(q),\dim D^3(q),\ldots\bigr),$$ where $D^j$ is the 
$j$-th power of the distribution $D$, i.e., 
$D^j=D^{j-1}+[D,D^{j-1}]$. A simple counting of the "number 
of parameters" in the considered equivalence problem shows 
that for $l=2$ the functional invariant should appear 
starting from $n=5$. It is well known that in the low 
dimensions $n=3$ or $4$ all generic germs of rank 2 
distributions are equivalent. (Darboux's theorem in the 
case $n=3$, small growth vector $(2,3)$ and Engel's theorem 
in the case $n=4$, small growth vector $(2,3,4)$, see, for 
example, \cite{bryantbook}, \cite{zhi}). 

In our previous paper \cite{zelduke}, using the notion of Jacobi curve 
of singular  extremals, introduced  in \cite{agrachev}, and general theory 
of unparametrized curves in Lagrange Grassmannian, developed in \cite{jac1} 
and \cite{jac2},
 we   
constructed  a special differential invariant, the {\it  
 fundamental form}
 of  generic rank 2 vector
distributions $D$ on $n$-dimensional manifold ($n\geq 5$). In the case $n=5$ this invariant can be realized as invariant homogeneous polynomial 
of degree $4$ on each plane $D(q)$ (in \cite{zelduke} we called this realization the {\it tangential fundamental form}). Our tangential fundamental form has the same algebraic nature as the {\it covariant binary biquadratic form}, constructed by E.Cartan in \cite{cartan}, using his ``reduction- prolongation'' procedure ( we call this form {\it Cartan's tensor}). 
In the present paper 
we prove that our  fundamental form coincides   
(up to constant factor $-35$) with Cartan's tensor.
Since 
these  two invariants were constructed in completely different ways, 
the comparison of 
them turned out to be not such an easy task.

The paper is organized as follows. In section 2, following 
\cite{zelduke}, we describe the construction of the 
fundamental form of (2,5)-distribution with the small 
growth vector $(2,3,5)$. In section 3, following chapter VI 
of the original paper \cite{cartan}, we briefly describe 
the main steps of construction of Cartan's tensor, 
rewriting all formulas that we need for the comparison of 
our and Cartan's invariants and formulate the main theorem 
of our paper (Theorem \ref{comptheor}).
  In section 4 
 for given frame, naturally adapted  to 
the distribution, we derive the formula for fundamental 
form in terms of its structural functions. 
This formula is important not only in the proof of our main 
theorem: in many cases it is more efficient from 
computational point of view than the method given in \cite 
{zelduke} ( see Theorem 2 there). Finally in section 5 we 
prove our main theorem. For this we just express our 
fundamental form in terms of structural functions of 
special adapted frame, distinguished by Cartan during his 
reduction process. 

In order to find the analog of Cartan's tensor in 
Riemannian geometry, one can apply the scheme $$ {\rm 
extremals}\,\, \rightarrow\,\, {\rm Jacobi\,\, 
curves}\,\,\rightarrow\,\, {\rm fundamental\,\, form\,\, 
of\,\, Jacobi\,\, curves}$$ to the (pseudo)-Riemannian 
metric and try to compare the obtained invariant with some 
classical invariants of Riemannian geometry. It was natural 
to expect that the obtained in this way invariant should be 
somehow related with Weyl conformal tensor of the metric. 
The exact relation was obtained by A. Agrachev and 
it will be presented in one of his forthcoming papers 
(compare with conjecture in 
\cite{bryant}, last paragraph of subsection 3.3 there).
 
{\it Acknowledgments.} I would like to thank professor Andrei 
 Agrachev for his constant attention to this work, stimulating discussions, and valuable advises. 
  I would also like to thank Boris Doubrov for using ChatGPT  in August 2026 to find a typo in the sign of one term in formula \eqref{FORMULA}. Namely, the term $\frac{1}{12}\vec h\circ\vec h\Bigl(b\, b_1+\vec h(b)-\frac{1}{3} a_3\Bigr)$ should have a negative sign. This can be easily traced by substituting the expression for $a_{21}$ from \eqref{a21f} into relation \eqref{auxA}. In the previous version, the sign of the term in \eqref{FORMULA} resulting from this substitution into the term $\frac{1}{12}a''_{21}$ of \eqref{auxA} was incorrect.


\section{Fundamental 
form of (2,5)-distribution } \indent 
\setcounter{equation}{0}
 
In this section
we describe the construction of the fundamental form of 
$(2,5)$-distribution with the small growth vector $(2,3,5)$. 
The presentation is close to \cite{zelduke}.

{\bf 2.1 Jacobi curve of abnormal extremals.} 
First for the distribution $D$ of the considered class  
one can distinguish special (unparametrized) curves in the 
cotangent bundle $T^*M$ of $M$. For this let 
$\pi:T^*M\mapsto M$ be the canonical projection. Let 
$\sigma$ be the standard symplectic structure on $T^*M$, 
namely, for any $\lambda\in T^*M$, $\lambda=(p,q)$, $q\in 
M$,$p\in T_q^*M$ let 
\begin{equation}
\label{stands}
\sigma(\lambda)(\cdot)=-d\, p\bigl(\pi_*(\cdot)\bigr)
\end{equation}
(here we prefer the sign ''-'' in the right handside , although usually one defines the standard symplectic form on $T^*M$ without this sign).
Denote 
by $(D^l)^{\perp}\subset T^*M$ the annihilator of the $l$th 
power $D^l$, namely 
\begin{equation}
\label{annihil} (D^l)^{\perp}= \{(q,p)\in T^*M:\,\, 
p\cdot v=0\,\,\forall v\in D^l(q)\}.
\end{equation}
The set $D^\perp$ is codimension $2$ submanifold of $T^*M$. 
Consider the restriction $\sigma|_{D^\perp}$ of the form 
$\sigma$ on $D^\perp$. It is not difficult to check that 
(see, for example \cite{zel}, section 2): the set of 
points, where the form $\sigma|_{D^\perp}$ is degenerated, 
coincides with $(D^2)^\perp$; the set 
$(D^2)^\perp\backslash(D^3)^\perp$ is codimension $1$ 
submanifold of $D^\perp$; for each $\lambda\in 
(D^2)^\perp\backslash(D^3)^\perp$ the kernel of 
$\sigma|_{D^\perp}(\lambda)$ is two-dimensional subspace of 
$T_\lambda D^\perp$, which is transversal to $T_\lambda 
(D^2)^\perp$. Hence $\forall\lambda\in 
(D^2)^\perp\backslash(D^3)^\perp$ we have $$ {\rm 
ker}\,\sigma|_{(D^2)^\perp}(\lambda)= {\rm 
ker}\,\sigma|_{D^\perp}(\lambda) \cap T_\lambda (D^2)^\perp 
$$ It implies that these kernels form line distribution in 
$(D^2)^\perp\backslash(D^3)^\perp$ and define a {\it 
characteristic 1-foliation} $Ab_D$ of 
$(D^2)^\perp\backslash(D^3)^\perp$. Leaves of this 
foliation will be called {\it characteristic curves} of 
distribution $D$. Actually these characteristic curves are 
so-called {\it abnormal extremals} of $D$. Projections of 
the characteristic curves to the base manifold $M$ will be 
called {\it abnormal trajectories} of $D$. Conversely, an 
abnormal extremal projected to the given abnormal 
trajectory will be called its {\it lift}.

\begin{remark}
\label{D30} {\rm 
 Note that in the 
considered case $(D^3)^\perp$ coincides with the zero 
section of $T^*M$.}$\Box$ \end{remark}

 Further, for a given 
segment $\gamma$ of characteristic curve one can construct 
a special (unparametrized) curve of Lagrangian subspaces, 
called Jacobi curve, in the appropriate symplectic space. 
For this for any $\lambda\in (D^2)^\perp$ denote by ${\cal 
J}(\lambda)$ the following subspace of 
$T_\lambda(D^2)^\perp$ 
\begin{equation}
\label{prejac} {\cal J}(\lambda)= \bigl(T_\lambda 
(T^*_{\pi(\lambda)}M)+ 
\ker\sigma|_{D^\perp}(\lambda)\bigr)\cap T_\lambda 
(D^2)^\perp.
\end{equation}
Here $T_\lambda (T^*_{\pi(\lambda)}M)$ is tangent to the 
fiber $T^*_{\pi(\lambda)}M$ at the point $\lambda$ (or 
vertical subspace of $T_\lambda(T^*M)$).
Actually in the considered case ${\cal J}$ is rank $4$ distribution 
on the manifold 
$(D^2)^\perp\backslash (D^3)^\perp$. 

 Let $O_\gamma$ be 
a neighborhood of $\gamma$ in $(D^2)^\perp$ such that 
\begin{equation}
\label{Ndef} N=O_\gamma /({Ab_D}|_{O_\gamma})
\end{equation}
 is a 
well-defined smooth manifold. The quotient manifold $N$ is 
a symplectic manifold endowed with a symplectic structure 
$\bar\sigma$ induced by $\sigma |_{(D^2)^\perp}$. Let $\phi 
:O_\gamma\to N$ be the canonical projection on the factor. 
It is easy to check that $\phi_*\bigl({\cal 
J}(\lambda)\bigr)$ is a Lagrangian subspace of the 
symplectic space $T_{\gamma}N$, $\forall\lambda\in\gamma$. 
 Let $L(T_\gamma 
N)$ be the Lagrangian Grassmannian of the symplectic space 
$T_\gamma N$, i.e., $$L(T_\gamma N)=\{\Lambda\subset 
T_\gamma N:\Lambda^\angle=\Lambda\},$$ where 
$\Lambda^\angle$ is the skew-symmetric complement of the 
subspace $\Lambda$, $$\Lambda^\angle=\{v\in T_\gamma 
N:\bar\sigma(v,\Lambda)=0\}.$$ {\it Jacobi curve of the 
characteristic curve (abnormal extremal)} $\gamma$ 
is the mapping 
\begin{equation}
\label{jacurve} \lambda\mapsto 
J_\gamma(\lambda)\stackrel{def}{=} \phi_*\Bigl({\cal 
J}(\lambda) 
\Bigr), \quad\lambda\in\gamma , 
\end{equation}
from $\gamma$ to $L(T_\gamma N)$.

 {\bf 2.2 Fundamental form of the curve in Lagrange 
Grassmannian} 
 Jacobi curves are invariants of the distribution $D$. They 
are unparametrized curves in the Lagrange Grassmannians.
 In 
\cite{jac1} for any curve of so-called constant {\it 
weight} in Lagrange Grassmannian we construct {\it the 
canonical projective structure} and a special degree 4  differential, 
{\it fundamental form}, which are invariants w.r.t. the action of 
linear Symplectic Group $GL(W)$ and reparametrization of the curve.  
Below we describe briefly  the construction of these invariants.
Actually it is more convenient to work with the curve in the set $G_m(W)$ be the set of all $m$-dimensional 
subspaces of $2m$-dimensional linear space  $W$ (i.e., in the Grassmannian of 
half-dimensional subspaces), where  General Linear Group acts.
 Since any curve 
of Lagrange subspaces w.r.t. some symplectic form on $W$ is 
obviously the curve in $G_m(W)$, all constructions below 
are valid for the curves in Lagrange Grassmannian.

For given $\Lambda\in G_m(W)$ denote by 
$\Lambda^\pitchfork$ the set of all $m$-dimensional 
subspaces of $W$ transversal to $\Lambda$, 
$$\Lambda^\pitchfork=\{\Gamma\in 
G_m(W):W=\Gamma\oplus\Lambda\}=\{\Gamma\in 
G_m(W):\Gamma\cap\Lambda=0\}$$ Fix some $\Delta\in 
\Lambda^\pitchfork$. Then for any subspace $\Gamma\in 
\Lambda^\pitchfork$ there exist unique linear mapping from 
$\Delta$ to $\Lambda$ with graph $\Gamma$. We denote this 
mapping by $\langle\Delta,\Gamma,\Lambda\rangle$. So, 
$$\Gamma=\{v+\langle\Delta,\Gamma,\Lambda\rangle 
v|v\in\Delta\}.$$ Choosing the bases in $\Delta$ and 
$\Lambda$ one can assign to any $\Gamma\in 
\Lambda^\pitchfork$ the matrix of the mapping
$\langle\Delta,\Gamma,\Lambda\rangle$ w.r.t these bases. In 
this way we define the coordinates on the set 
$\Lambda^\pitchfork$. 

\begin{remark}
\label{lagr} {\rm Assume that $W$ is provided with some 
symplectic form $\bar\sigma$ and $\Delta, \Lambda$ are 
Lagrange subspaces w.r.t. $\bar\sigma$. Then the map 
$v\mapsto \bar\sigma(v, \cdot)$, $v\in \Delta$, defines the 
canonical isomorphism between $\Delta$ and $\Lambda^*$. It 
is easy to see that $\Gamma$ is Lagrange subspace iff the 
mapping $\langle\Delta,\Gamma,\Lambda\rangle$, considered 
as the mapping from $\Lambda^*$ to $\Lambda$, is 
self-adjoint.} $\Box$ \end{remark} 

Let $\Lambda(t)$ be a smooth curve in $G_m(W)$ defined on 
some interval $I\subset\mathbb{R}$. We are looking for 
invariants of $\Lambda(t)$ by the action of $GL(W)$. 
We say that the curve $\Lambda(\cdot)$ is {\it ample at} 
$\tau$ if $\exists s>0$ such that for any representative 
$\Lambda^s_\tau(\cdot)$ of the $s$-jet of $\Lambda(\cdot)$ 
at $\tau$, $\exists t$ such that 
$\Lambda^s_\tau(t)\cap\Lambda(\tau)=0$. The curve 
$\Lambda(\cdot)$ is called {\it ample} if it is ample at 
any point. This is an intrinsic definition of an ample 
curve. In coordinates this definition takes the following 
form: if in some coordinates the curve $\Lambda(\cdot)$ is 
a curve of matrices $t\mapsto S_t$, then $\Lambda(\cdot)$ 
is ample at $\tau$ if and only if the function $t\mapsto 
det(S_t-S_\tau)$ has a root of {\sl finite order} at 
$\tau$. 

\begin{defin}
The order of zero of the function $t\mapsto 
det(S_t-S_\tau)$ at $\tau$, where $S_t$ is a coordinate 
representation of the curve $\Lambda(\cdot)$, is called a 
weight of the curve $\Lambda(\cdot)$ at $\tau$. 
\end{defin}  
It is clear that the weight of $\Lambda(\tau)$ is integral 
valued upper semicontinuous functions of $\tau$. Therefore 
it is locally constant on the open dense subset of the 
interval of definition $I$. 

Now suppose that the curve has the constant weight $k$ on 
some subinterval $I_1\subset I$. It implies that for all 
two parameters $t_0$,$t_1$ in $I_1$ sufficiently such that 
$t_0\neq t_1$ , one has 
 $$
 \Lambda(t_0)\cap\Lambda(t_1)=0.
 $$
 Hence for such $t_0$, $t_1$ the following linear mappings
 \begin{eqnarray} ~&~
\frac{d}{d s}\langle\Lambda(t_0),\Lambda(s),\Lambda(t_1)
\rangle\Bigl|_{s=t_0}\Bigr 
.:\Lambda(t_0)\mapsto\Lambda(t_1),\label{op1}\\ 
 ~&~\frac{d}{d s}\langle\Lambda(t_1),\Lambda(s),\Lambda(t_0)
 \rangle\Bigl|_{s=t_1}\Bigr 
 .:\Lambda(t_1)\mapsto\Lambda(t_0)\label{op2}
 \end{eqnarray}
 are well defined. Taking composition of mapping (\ref{op2}) 
 with mapping  (\ref{op1}) we obtain the operator 
 from the subspace $\Lambda(t_0)$ to itself.
 
 \begin{remark}
 \label{crossrem}
 {\rm  This operator  is 
 actually the infinitesimal cross-ratio of two points
 $\Lambda(t_i)$, $i=0,1$, together with two tangent vectors 
 $\dot\Lambda(t_i)$, $i=0,1$, at these points in $G_m(W)$ 
 (see \cite{jac1} for the details).}
 $\Box$
 \end{remark}

\begin{theor} (see \cite{jac1}, Lemma 4.2) If the curve has the constant weight $k$ on 
some subinterval $I_1\subset I$, then the following 
asymptotic holds 
 \begin{eqnarray}
  ~&~{\rm tr}\Bigl( \frac{d}{ds}\langle\Lambda(t_1),\Lambda(s),\Lambda(t_0)
 \rangle\Bigl|_{s=t_1}\Bigr.\circ
  \frac{d}{ds}\langle\Lambda(t_0),
 \Lambda(s),\Lambda(t_1)
 \rangle\Bigl|_{s=t_0}\Bigr.\Bigr)=\nonumber\\
 ~~\label{mainass}\\
 ~&~-\frac {k}{(t_0-t_1)^2}-g_\Lambda(t_0,t_1),\nonumber
 \end{eqnarray}
 where $g_{_{\Lambda}}(t,\tau)$ 
  is a smooth function in the neighborhood of diagonal 
  $\{(t,t)| t\in I_1\}$. 
\end{theor}

The function $g_{_{\Lambda}}(t,\tau)$ is "generating" 
function for invariants of the parametrized curve by the 
action of $GL(2m)$. The first coming invariant of the 
parametrized curve, {\it the generalized Ricci curvature}, is 
just $g_{_{\Lambda}}(t,t)$, 
 the value of $g_{_{\Lambda}}$  at the diagonal. 

In order to obtain invariants for unparametrized curves 
(i.e., for one-dimensional submanifold of $G_m(W)$) we use 
a simple reparametrization rule for a function 
$g_{_{\Lambda}}$. Indeed, let 
$\vf:\mathbb{R}\mapsto\mathbb{R}$ be a smooth monotonic 
function. It follows directly from (\ref{mainass}) that 
\begin{equation}\label{chainrul} g_{_{\Lambda\circ 
\vf}}(t_0,t_1)=\dot\vf(t_0)\dot\vf(t_1)g_{_{\Lambda}} 
(\vf(t_0),\vf(t_1))+k\left(\frac{\dot\vf(t_0)\dot\vf(t_1)} 
{(\vf(t_0)-\vf(t_1))^2}-\frac{1}{(t_0-t_1)^2}\right). 
\end{equation}
In particular, putting $t_0=t_1=t$, one obtains the 
following reparametrization rule for the generalized Ricci 
curvature 
\begin{equation}
 \label{rhorep}
 g_{_{\Lambda\circ 
\vf}}(t,t)=\dot\vf(t)^2g_{_{\Lambda}}(\vf(t),\vf(t))+ 
\frac{k}{3}\mathbb{S}(\vf),
\end{equation}
 where $\mathbb{S}(\vf)$ is a 
Schwarzian derivative of $\vf$, \begin{equation} 
\label{sch} \mathbb S(\varphi)= 
\frac{1}{2}\frac{\varphi^{(3)}}{\varphi'}- \frac{3}{4}\Bigl 
(\frac{\varphi''}{\varphi'}\Bigr)^2=\frac 
{d}{dt}\Bigl(\frac {\varphi''}{2\,\varphi'}\Bigr) 
-\Bigl(\frac{\varphi''}{2\,\varphi'}\Bigr)^2.
\end{equation}

From (\ref{rhorep}) it follows that the class of local 
parametrizations, in which the generalized Ricci curvature 
is identically equal to zero, defines a {\it canonical 
projective structure} on the curve (i.e., any two 
parametrizations from this class are transformed one to 
another by M{\"obius transformation). This parametrizations 
are called projective. From (\ref{chainrul}) it follows 
that if $t$ and $\tau$ are two projective parametrizations 
on the curve $\Lambda(\cdot)$, 
$\tau=\vf(t)=\frac{at+b}{ct+d}$, and $g_\Lambda$ is 
generating function of $\Lambda(\cdot)$ w.r.t. parameter 
$\tau$ then 
\begin{equation}
\label{prefund} {\partial^2 
g_{_{\Lambda\circ\vf}}\over\partial 
t_1^2}(t_0,t_1)\Bigr|_{t_0=t_1=t}\Bigl.={\partial^2 
g_{_{\Lambda}}\over\partial 
\tau_1^2}(\tau_0,\tau_1)\Bigr|_{\tau_0=\tau_1=\vf(t)}\Bigl. 
(\vf'(t))^4, 
\end{equation}
which implies that the following degree four differential 
$$ {\cal A}=\frac{1}{2} {\partial^2 
g_{_{\Lambda}}\over\partial 
\tau_1^2}(\tau_0,\tau_1)\Bigr|_{\tau_0=\tau_1=\tau}\Bigl. 
(d\tau)^4$$ on the curve $\Lambda(\cdot)$ does not depend 
on the choice of the projective parametrization (by degree 
four differential on the curve we mean the following: for 
any point of the curve a degree $4$ homogeneous function is 
given on the tangent line to the curve at this point). This 
degree four differential is called a {\it fundamental form} 
of the curve. 
If $t$ is an arbitrary (not necessarily projective) 
parametrization on the curve $\Lambda(\cdot)$, then the 
fundamental form in this parametrization has to be of the 
form $A(t)(dt)^4$, where $A(t)$ is a smooth function, 
called the {\it density} of the fundamental form w.r.t. 
parametrization $t$. 

Now we suppose that the linear space $W$ is provided with 
some symplectic form $\bar \sigma$ and we restrict 
ourselves to the curves in the corresponding Lagrange 
Grassmannian $L(W)$. The tangent space $T_\Lambda L(W)$ to 
any $\Lambda\in L(W)$ can be identified with the space of 
quadratic forms $Q(\Lambda)$ on the linear space $\Lambda$. 
Indeed, take a curve $\Lambda(t)\in G_k(W)$ with 
$\Lambda(0)=\Lambda$. Given some vector $l\in\Lambda$, take 
a curve $l(\cdot)$ in $W$ such that 
 $l(t)\in \Lambda (t)$ for all sufficiently small $t$ 
 and $l(0)=l$.
 It is 
easy to see that the quadratic form $l\mapsto 
\bar\sigma\bigl(l'(0),l\bigr)$ depends only on 
$\frac{d}{dt}\Lambda(0)$. In this way we identify 
$\frac{d}{dt}\Lambda(0)\in T_\Lambda G_k(W)$ with some 
element of $Q(\Lambda)$ (a simple counting of dimension 
shows that these correspondence between $T_\Lambda L(W)$ 
and 
 $Q(\Lambda)$ is a bijection).} 
The rank of the curve  $\Lambda(t)\subset L(W)$ is by definition the rank of 
its velocity $\frac{d}{dt}\Lambda(t)|_{t=\tau}$ at $\tau$, considered as the quadratic form. The curve $\Lambda(t)$ in $L(W)$ is 
called
{\it nondecreasing (nonincreasing)} if its velocities $\frac{d}{dt}\Lambda(t)$ at any point are nonnegative (nonpositive) definite quadratic forms.


{\bf 2.3 Construction of  fundamental form of (2,5)- distributions} 
  Note that Jacobi curve $J_\gamma$ of characteristic curve $\gamma$ 
 of distribution $D$
 defined by (\ref{jacurve}) is not 
ample, because all subspaces $J_\gamma(\lambda)$ have a 
common line. Indeed, let $\delta_a:T^*M\mapsto T^*M$ be the 
homothety by $a\neq 0$ in the fibers, namely, 

\begin{equation}
\label{homoth}
 \delta_a(p,q)=(ap,q), \quad q\in M,\,\,p\in 
T^*M.
\end{equation}
 Denote by $\vec e(\lambda)$ the following vector 
field called Euler field \begin{equation} \label{Euler} 
\vec e(\lambda)= \frac{\partial}{\partial a}\delta
_a(\lambda)\Bigl|_{a=1}\Bigr. 
\end{equation} 
\begin{remark}
\label{homrem}{\rm Obviously, if $\gamma$ is characteristic 
curve of $D$, 
then also $\delta_a(\gamma)$ is.} $\Box$ 
\end{remark}

It implies that the vectors $\phi_*(\vec 
e(\lambda))$ coincide for all $\lambda\in \gamma$, so the 
line 
\begin{equation}
\label{ega} E_\gamma\stackrel{def}{=}\{\mathbb R 
\phi_*(\vec e(\lambda))\}
\end{equation}
 is common for all subspaces 
$J_\gamma(\lambda)$, $\lambda\in \gamma$ (here, as in 
Introduction, $\phi 
:O_\gamma\to N$ is the canonical projection on the factor 
$N=O_\gamma /({Ab_D}|_{O_\gamma})$, where $O_\gamma$ is 
sufficiently small tubular neighborhood of $\gamma$ in 
$(D^2)^\perp$). 

Therefore it is natural to make an appropriate 
factorization by this common line $E_\gamma$. Namely, by 
above all subspaces $J_\gamma(\lambda)$ belong to 
skew-symmetric complement $E_\gamma^\angle$ of $E_\gamma$ 
in $T_\gamma N$. Denote by $p:T_\gamma N\mapsto T_\gamma 
N/E_\gamma$ the canonical projection on the factor-space. 
The mapping 
\begin{equation}
\label{redjac} \lambda\mapsto \widetilde 
J_\gamma(\lambda)\stackrel{def}{=}p(J_\gamma(\lambda)), 
\,\,\lambda\in\gamma 
\end{equation}
from $\gamma$ to $L(E_\gamma^\angle/E_\gamma)$ is called 
{\it reduced Jacobi curve} of characteristic curve 
$\gamma$. 
Note that in the considered case 
\begin{equation}
\label{dimjac}
{\rm dim }\,E_\gamma^\angle/E_\gamma=4
\end{equation}

To clarify the notion of the reduced Jacobi curve note that 
for any $\lambda\in\gamma$ one can make the following 
identification 
\begin{equation}
\label{tanid} T_\gamma N\sim T_\lambda(D^2)^\perp/{\rm 
ker}\, \sigma|_{(D^2)^\perp}(\lambda). 
\end{equation}
Take on $O_\gamma$ any vector field $H$ 
 tangent to characteristic $1$-foliation $Ab_D$ and without stationary points,
i.e., $H(\lambda)\in {\rm ker} 
\sigma|_{(D^2)^\perp}(\lambda)$, $H(\lambda)\neq 0$ for all 
$\lambda\in O_\gamma$. Then it is not hard to see that 
under identification (\ref{tanid}) one has 
\begin{equation}
\label{ident} \widetilde J_\gamma(e^{t H}\lambda)=(e^{-t 
H})_ 
*\Bigl( {\cal J}(e^{t 
H}\lambda) 
\Bigr) /{\rm span}( {\rm 
ker}\,\sigma|_{(D^2)^\perp}(\lambda),\vec e(\lambda) ) 
\end{equation} 
where $e^{t H}$ is the flow generated by the vector field $ 
H$, and the subspaces $\widetilde J_\gamma(e^{t H}\lambda)$ 
live in the symplectic space $W_\lambda$, defined as 
follows
\begin{equation}
\label{wl} W_\lambda=\Bigl(\vec e(\lambda)^\angle\cap 
T_\lambda (D^2)^\perp\Bigr)/{\rm span} \bigr({\rm ker}\, 
\sigma|_{(D^2)^\perp}(\lambda) , \vec e(\lambda)\bigl).
\end{equation} 

The following proposition follows directly from Propositions 2.2, 2.6 of 
\cite{zelduke} and relation (\ref{dimjac}).
\begin{prop}
\label{ranklem} 
 The reduced Jacobi curve of characteristic curve of $(2,5)$-distribution 
with the small growth vector $(2,3,5)$ is rank 1 
nondecreasing curve of the constant weight $4$ in Lagrange 
Grassmannian of $4$-dimensional linear symplectic space. 
\end{prop}

In order to construct fundamental form we have to introduce some notations.
Let  $X_1$, $X_2$ be two vector fields, constituting
the basis of distribution $D$,
 i.e.,
\begin{equation}
\label{X12} D( q)={\rm span}(X_1(q), X_2(q))\,\,\,\,\forall 
q\in 
M.
\end{equation}
Since our study is local, we can always suppose that such 
basis exists, restricting ourselves, if necessary, on some 
coordinate neighborhood instead of whole $M$. 
Given the basis  $X_1$, $X_2$  one can construct a special 
vector field $\vec h_{_{X_1,X_2}}$ tangent to characteristic $1$-foliation 
$Ab_D$.
For this suppose that 
\begin{eqnarray}
&~&X_3=[X_1,X_2]\quad {\rm mod}\, D,\,\,\, 
X_4=[X_1,[X_1,X_2]]=[X_1,X_3]\quad {\rm mod}\, 
D^2,\nonumber 
\\ &~&~\label{x345} \\ 
&~&X_5=[X_2,[X_1,X_2]]=[X_2,X_3]\quad{\rm mod}\, D^2 
\nonumber 
\end{eqnarray}
\begin{defin}
\label{adapt} The tuple $\{X_i\}_{i=1}^5$, satisfying 
(\ref{X12}) and (\ref{x345}) is called adapted frame of 
the distribution $D$. If instead of (\ref{x345}) one has 
\begin{equation}
\label{x345s} X_3=[X_1,X_2],\,\,\, 
X_4=[X_1,[X_1,X_2]]=[X_1,X_3],\,\,\,\, 
X_5=[X_2,[X_2,X_1]]=[X_3,X_2], 
\end{equation}
the frame $\{X_i\}_{i=1}^5$ will be called {\it 
strongly adapted} to $D$.
\end{defin}
 
Let us introduce ``quasi-impulses'' $u_i:T^*M\mapsto\mathbb 
R$, $1\leq i\leq 5$, 
\begin{equation}
\label{quasi25} u_i(\lambda)=p\cdot X_i(q),\,\,\lambda=(p,q),\,\,
q\in M,\,\, p\in T_q^* M
\end{equation}
For given function $G:T^*M\mapsto \mathbb R$ denote by 
$\vec G$ the corresponding Hamiltonian vector field defined 
by the relation $\sigma(\vec G,\cdot)=d\,G(\cdot)$.
Then it is easy to show (see, for example \cite{zel})
that
\begin{equation}
\label{ker25} \ker\sigma\Bigr|_{D^\perp}\Bigl.(\lambda)= 
{\rm span}(\vec
u_1(\lambda),\vec u_2(\lambda)),\quad \forall \lambda\in D^\perp, 
\end{equation} 
\begin{equation}
\label{foli25}
\ker\sigma\Bigr|_{(D^2)^\perp}\Bigl.(\lambda)= \mathbb{R}
\Bigl((u_4 \vec{u}_2-u_5\vec{u}_1)(\lambda)\Bigr),\quad 
\forall \lambda\in (D^2)^\perp\backslash (D^3)^\perp
\end{equation}
The last relation implies that 
the following vector field 
\begin{equation}
\label{ham25} \vec h_{_{X_1,X_2}}=u_4 \vec u_2- u_5 \vec 
u_1 
\end{equation}
is tangent to the characteristic $1$-foliation.

Now we are ready to construct the fundamental form.
For any $\lambda\in (D^2)^\perp\backslash (D^3)^\perp$ take 
characteristic curve $\gamma$, passing through $\lambda$. 
Let ${\cal A}_\lambda$ be the fundamental form of the 
reduced Jacobi curve $\widetilde J_\gamma(\lambda)$ of 
$\gamma$ at $\lambda$ (by Proposition \ref{ranklem} ${\cal A}_\lambda$ is well defined). By construction ${\cal A}_\lambda$ 
is degree $4$ homogeneous function on the tangent line to 
$\gamma$ at $\lambda$. Given a local basis
$(X_1, X_2)$ of distribution $D$ let
\begin{equation}
\label{densX12} A_{_{X_1,X_2}}(\lambda)= {\cal 
A}_\lambda(\vec h_{_{X_1,X_2}}(\lambda)) \end{equation} In 
this way 
 to any (local) basis 
$(X_1, X_2)$ of distribution $D$ we assign the function $ 
A_{_{X_1,X_2}}$ on $(D^2)^\perp\backslash (D^3)^\perp$. 

\begin{remark}
\label{densrem} {\rm If we consider parametrization 
$t\mapsto \widetilde J_\gamma(e^{t\vec 
h_{_{X_1,X_2}}}\lambda)$ of the reduced Jacobi curve of 
$\gamma$, then $A_{_{X_1,X_2}}(\lambda)$ is the density of 
fundamental form of this curve w.r.t. parametrization $t$ 
at $t=0$.} $\Box$ 
\end{remark}

Take another basis $\tilde X_1,\tilde X_2$ of the 
distribution $D$. It can be shown (see \cite {zelduke}, 
beginning of subsection 2.3) that for any $q\in M$ the 
restriction of the corresponding function $ A_{_{\tilde 
X_1,\tilde X_2}}$ to the fiber $(D^2)^\perp(q)\backslash 
(D^3)^\perp(q)$ ($=(D^2)^\perp(q)\backslash \{0\}$ by 
Remark \ref{D30}) is equal to the restriction of $ 
A_{_{X_1,X_2}}$ to $(D^2)^\perp(q)\backslash \{0\}$, 
multiplied on some positive constant (the square of the
determinant of transition matrix from the basis $(X_1,X_2)$ 
to the basis $(\tilde X_1,\tilde X_2)$). The following 
Proposition is direct consequence of the last assertion, 
Proposition 2.8 and Theorem 3 from \cite{zelduke}: 

\begin{prop}
\label{fundform} If $D$ is $(2,5)$ distribution on $M$ with 
the small growth vector $(2,3,5)$ at any point, then the 
restriction of $A_{_ {X_1,X_2}}$ to the fiber 
$(D^2)^\perp(q)$ is well defined degree 4 homogeneous 
polynomial, up to multiplication on positive constant. 
\end{prop}
 
\begin{defin}
\label{fundefin}
The restriction of $A_{_ {X_1,X_2}}$ to $(D^2)^\perp(q)$ will be called
 {\it fundamental form of $(2,5)$-distribution D at the point $q$}.
\end{defin}
From now on we will 
write $\vec h$ instead of $\vec h_{_{X_1,X_2}}$ and $A$ 
instead of $A_{_{X_1,X_2}}$ without special mentioning. 

In the case $n=5$ and small growth vector $(2,3,5)$ one can 
look at the fundamental form of the distribution $D$ from 
the different point of view. In this case (in contrast to 
generic $(2,n)$-distributions with $n>5$) there is only one 
abnormal trajectory starting at given point $q\in M$ in 
given direction (tangent to $D(q)$). All lifts of this 
abnormal trajectory can be obtained one from another by 
homothety. So they have the same, up to symplectic 
transformation, Jacobi curve. It means that one can 
consider Jacobi curve and fundamental form of this curve on 
abnormal trajectory instead of abnormal extremal. 
Therefore, to any $q\in M$ one can assign a homogeneous 
degree $4$ polynomial $\AA_q$ on the plane $D(q)$ in 
the following way: 
\begin{equation}
\label{tangfund} \AA_q(v)\stackrel{def}{=}{\cal 
A}_\lambda(H) 
\end{equation} for any $v\in D(q)$, where 
\begin{equation}
\label{condlh} \pi(\lambda)=q,\quad\pi_*H=v,\quad H\in {\rm 
ker}\,\sigma|_{(D^2)^\perp}(\lambda). \end{equation} and 
the righthand side of (\ref{tangfund}) is the same for any 
choice of $\lambda$ and $H$, satisfying (\ref{condlh}). 

$\AA_q$ will be called {\it tangential fundamental form} of 
the distribution $D$ at the point $q$. We stress that the 
tangential fundamental form is the well defined degree $4$ 
homogeneous polynomial on $D(q)$ and not the class of 
polynomials defined up to multiplication on a positive 
constant. 

In \cite{zelduke} in the case $n=5$ we obtained the 
explicit formulas for calculation of the (tangential) 
fundamental form (see Theorem 2 there), but they are not 
sufficient for our purposes here.

\section{Cartan's tensor}
\indent
\setcounter{equation}{0} 

In the present section, following chapter VI of the original 
paper \cite{cartan}, we briefly describe the main steps of 
construction of Cartan's tensor, rewriting all formulas 
that we need for the comparison of our and Cartan's 
invariants. We will use the language of Cartan and his 
notations, referring sometimes by remarks to modern 
terminology of $G$-structures (see 
\cite{gard},\cite{mont},\cite{stern}). 
In order to simplify the formulas,
 we will omit in the sequel the sign $\wedge$ 
 for standard operation with differential forms.
 
Let $\omega_1$, $\omega_2$, $\omega_3$, $\omega_4$, and 
$\omega_5$ be coframe on $M$ (i.e., $\{\omega 
_i(q)\}_{i=1}^5$ constitute the basis of $T_q^*M$ for any 
$q\in M$) such that the rank $2$ distribution $D$ is 
annihilator of the first three elements of this coframe, 
namely, 

\begin{equation}
\label{annih}
 D(q)=\{v\in T_qM;\,\, 
\omega_1(v)=\omega_2(v)=\omega_3(v)=0\},\quad \forall q\in 
M 
\end{equation}

Obviously the set of all coframes satisfying (\ref{annih}) 
is 19-parametric family. Among all these coframes Cartan 
distinguishes special 
coframes 
satisfying the following structural equations (formula 
(VI.5) in \cite{cartan}): 

\begin{eqnarray}
&~&d\,\omega_1=\omega_1(2\overline\omega_1+\overline\omega_4)+ 
\omega_2\overline\omega_2+\omega_3\omega_4\nonumber\\ 
&~&d\,\omega_2=\omega_1\overline\omega_3+\omega_2 ( 
\overline\omega_1+2\overline\omega_4)+\omega_3\omega_5\nonumber\\ 
&~&d\,\omega_3=\omega_1\overline\omega_5+ 
\omega_2\overline\omega_6+\omega_3(\overline\omega_1 
+\overline\omega_4)+\omega_4\omega_5\label{costruct}\\ 
&~&d\,\omega_4=\omega_1\overline\omega_7+ 
\frac{4}{3}\omega_3\overline\omega_6+ 
\omega_4\overline\omega_1+\omega_5\overline\omega_2\nonumber\\ 
&~&d\,\omega_5=\omega_2\overline\omega_7- 
\frac{4}{3}\omega_3\overline\omega_5+ 
\omega_4\overline\omega_3+\omega_5\overline\omega_4,\nonumber 
\end{eqnarray}
where $\overline \omega_j$, $1\leq j\leq 7$, are new 
$1$-forms. It turns out that the set of all coframes
satisfying (\ref{costruct}) is $7$-parametric family. This 
family defines $12$-dimensional bundle over $M$ , which 
will be denoted by ${\cal K}(M)$: instead of considering a 
family of coframes $\{\om_i\}_{i=1}^5$ on $M$ one can 
consider the $5$-tuple of $1$-forms on ${\cal K}(M)$; one 
can think of the forms $\overline \omega_j$, $1\leq j\leq 
7$, from (\ref{costruct}) and of the equation 
(\ref{costruct}) itself as defined on ${\cal K}(M)$. 
\begin{remark}{\rm
In the modern terminology ${\cal K}(M)$ is 
$G$-structure over $M$ such that the Lie algebra of its 
structure group is the algebra of the following matrices 
\begin{equation}
\label{g7} \left(\begin{array}{ccccc} 2a_1+a_4& a_2&0&0&0\\ 
a_3&a_1+2a_4&0&0&0\\ a_5&a_6&a_1+a_4&0&0\\ 
a_7&0&\frac{4}{3}a_6&a_1&a_2\\ 0&a_7&-\frac{4}{3} 
a_5&a_3&a_4 
\end{array}
\right)
\end{equation} 
(compare this with (\ref{costruct})).One says that the 
$G$-structure ${\cal K}(M)$ is obtained by reduction 
procedure from the bundle of all coframes satisfying 
(\ref{annih}). It turns out that the further reduction is 
impossible.} $\Box$ 
\end{remark}
Note also that the tuple of the forms
$\Bigl(\{\om_i\}_{i=1}^5,\{\bom_j\}_{j=1}^7\Bigr)$ is 
coframe on ${\cal K}(M)$. Further, it turns out that the 
tuple of the forms $\overline \omega_j$, $1\leq j\leq 7$, 
is defined by structure equation (\ref{costruct}) up to the 
following transformations (the formula after the formula 
(VI.7) in \cite{cartan}): 
\begin{equation}
\label{nu} 
\begin{array}{c}
\Bigl(\overline \omega_1,\,\overline \omega_2,\,\overline \omega_3,\,\overline \omega_4,\,\overline \omega_5,\,\overline \omega_6,\,\overline \omega_7\Bigr)\rightarrow
\Bigl(\overline\omega_1+\nu_1\omega_1,\,\overline\omega_2+\nu_2\omega_1,\\ 
\overline\omega_3+\nu_1\omega_2, \,\overline\omega_4+\nu_2\omega_2, 
\overline\omega_5+\nu_1\omega_3,\,\overline\omega_6+\nu_2\omega_3,\,
\overline\omega_7+\nu_1\omega_4+\nu_2\overline\omega_5\Bigr),
\end{array}
\end{equation}  
where parameters $\nu_1$ and $\nu_2$ are arbitrary. 
Replacing the tuple $\overline \omega_j$, $1\leq j\leq 7$, 
by the righthand side of (\ref{nu}), one has $2$-parametric 
family of coframes 
$\Bigl(\{\om_i\}_{i=1}^5,\{\bom_j\}_{j=1}^7\Bigr)$ on 
${\cal K}(M)$. This family defines $14$-dimensional bundle 
over ${\cal K}(M)$, which will be denoted by ${\cal 
K}_1(M)$. 
\begin{remark}{\rm
One says that the bundle ${\cal K}_1(M)$ is prolongation of 
${\cal K}(M)$.}$\Box$ 
\end{remark}
Further Cartan expresses seven exterior derivatives 
$d\,\bom$, considered as the forms on ${\cal K}_1(M)$, by 
the $1$-forms $\{\om_i\}_{=1}^5$, 
$\{\bom_j\}_{j=1}^7$
and two new forms 
$\chi_1$ and $\chi_2$ on ${\cal K}_1(M)$ such that these 
two new forms are defined uniquely by these expressions(see 
formula (VI.8) in \cite{cartan}). Therefore the tuple 
\begin{equation}
\label{cartcan}
\Bigl(\{\om_i\}_{i=1}^5,\{\bom_j\}_{j=1}^7,\chi_1,\chi_2\Bigr)
\end{equation}
is canonical coframe on ${\cal K}_1(M)$. 

\begin{remark}
{\rm The pair ${\cal K}_1(M)$ and coframe (\ref{cartcan}) 
can be considered as the trivial bundle ${\cal K}_2(M)$ 
over ${\cal K}_1(M)$, which is actually one more 
prolongation. So, by reduction from the bundle of coframes 
satisfying (\ref{annih}) to the bundle ${\cal K}(M)$ and by 
two successive prolongations (from ${\cal K}(M)$ to ${\cal 
K}_1(M)$ and then to ${\cal K}_2(M)$) one can arrive to the 
canonical coframe, which essentially solve the problem of 
equivalence of the considered class of distributions.} 
$\Box$ 
\end{remark}
In the sequel we will need the following consequences of 
the formula (VI.8) in \cite{cartan}: 
\begin{eqnarray}
d\,(\overline \omega_1-\bom_4)&=&2\bom_3\bom_2 -
\omega_4\overline \omega_5 +\omega_5\overline 
\omega_6+\omega_1 \chi_1-\om_2\chi_2+3 B_2 
\omega_1\omega_3+ \nonumber\\ 
~&~& 3 B_3 \omega_2\omega_3+3 A_2 \om_1\om_4+3 
A_3\om_1\om_5+ 3 A_3\om_2\om_4+ 3 A_4\om_2\om_5, 
\label{bom1}\\ 
 d\,\bom_2&=&\bom_2(\bom_1-\bom_4)-\om_4\bom_6 
+\om_1\chi_2 + B_4\om_2\om_3+A_4\om_2\om_4+ 
A_5\om_2\om_5,\label{bom2}\\
d\,\bom_3&=&\bom_3(\bom_4-\bom_1)-\om_5\bom_5 +\om_2\chi_1 
- B_1\om_1\om_3 -A_1\om_1\om_4- A_2\om_1\om_5,\label{bom3} 
\end{eqnarray}
Here $A$'s and $B$'s
are functions on ${\cal K}_1(M)$.

Cartan considers the following expression 
\begin{equation}
\label{cartf} {\cal F}=A_1 \om_4^4+4 A_2\om_4^3\om_5 +6 
A_3\om_4^2\om_5^2+4 A_4\om_4\om_5^3+A_5\om_5^4. 
\end{equation}
Let ${\cal P}:{\cal K}_1(M)\mapsto M$ be the canonical 
projection. Note that for given $Q\in {\cal K}_1(M)$ the 
form ${\cal F}$ at $Q$ can be considered as degree $4$
homogeneous polynomial on $T_{{\cal P}(Q)}M$. It turns out 
that the form ${\cal F}$ calculated at different points of 
the same fiber ${\cal P}^{-1}(q)$, $q\in M$, gives the same 
polynomial (modulo $\om_1(q)$, $\om_2(q)$, $\om_3(q)$) on 
$T_q M$ or, equivalently, the same polynomial on the plane 
$D(q)$ (recall that by construction $D(q)$ is annihilator 
of ${\rm span} (\om_1(q), \om_2(q), 
\om_3(q)$). 
Briefly speaking, ${\cal F}$ restricted on $D$ is covariant 
symmetric tensor of order $4$. This tensor is called {\it 
Cartan's tensor} of distribution $D$. We will denote by 
${\cal F}_q$ Cartan's tensor at the point $q$. 

So for any $q\in M$ our tangential fundamental form $\AA_q$ 
and Cartan's tensor ${\cal F}_q$ are both invariantly 
defined degree 4 homogeneous polynomials on the plane 
$D(q)$. The following theorem is the main results of the 
present paper: 

\begin{theor}
\label{comptheor} For any $q\in M$ Cartan's tensor ${\cal 
F}_q$ and the tangential fundamental form $\AA_q$ are 
connected by the following identity 
\begin{equation}
\label{FAA} {\cal F}_q= -35 \AA_q. \end{equation} 
\end{theor} 

As a conclusion of this theorem we obtain 
that  geometric reason for the existence of 
Cartan's tensor is the existence of the special degree 4 
differential on curves in Grassmannian of half-dimensional 
subspaces, which is constructed with the help of the notion 
of the cross-ratio of four point in this Grassmanian (or of 
the infinitesimal cross-ratio of two points in this 
Grassmannian together with two tangent vectors at these 
points). 
\section{Fundamental form for $n=5$ in terms of the structural functions 
of some adapted frame} \indent \setcounter{equation}{0}

In the present section we make preparations to prove 
Theorem \ref{comptheor}. Namely, for given adapted frame to 
the distribution we derive the formula for fundamental form 
in terms of its structural functions. 
This formula is important not only in the proof of this 
theorem: in many cases it is more efficient from 
computational point of view than the method given in \cite 
{zelduke} ( see Theorem 2 there).

First  
we need more facts from the theory of curves in 
Grassmannian $G_m W)$ of half-dimensional subspaces (here 
${\rm dim}\, W=2m$) and in Lagrange Grassmannian $L(W)$ 
w.r.t. to some symplectic form on $W$, developed in 
{\cite{jac1}, \cite{jac2}). Below we present all necessary 
facts from the mentioned papers together with several new 
useful arguments. 

Fix some $\Lambda\in G_m(W)$. As before , let 
$\Lambda^\pitchfork$ be the set of all $m$-dimensional 
subspaces of $W$ transversal to $\Lambda$. Note that any 
$\Delta\in \Lambda^\pitchfork$ can be canonically 
identified with $W/\Lambda$. Keeping in mind this 
identification and taking another subspace $\Gamma\in 
\Lambda^\pitchfork$ one can define the operation of 
subtraction $\Gamma-\Delta$ as follows 
$$\Gamma-\Delta\stackrel{def}{=}\langle\Delta,\Gamma,\Lambda\rangle\in 
{\rm Hom}\,(W/\Lambda,\Lambda).$$ It is clear that the set 
$\Lambda^\pitchfork$ provided with this operation can be 
considered as the affine space over the linear space ${\rm 
Hom}\, (W/\Lambda,\Lambda)$.

Consider now some ample curve $\Lambda(\cdot)$ in $G_m(W)$. 
Fix some parameter $\tau$. By assumptions 
$\Lambda(t)\in\Lambda(\tau)^\pitchfork$ for all $t$ from a 
punctured neighborhood of $\tau$. We obtain the curve 
$t\mapsto\Lambda(t)\in\Lambda(\tau)^\pitchfork$ in the 
affine space $\Lambda(\tau)^\pitchfork$ with the pole at 
$\tau$. We denote by $\Lambda_\tau(t)$ the identical 
embedding of $\Lambda(t)$ in the affine space 
$\Lambda(\tau)^\pitchfork$. First note that the velocity 
${\partial\over\partial t}\Lambda_\tau(t)$ is well defined 
element of ${\rm Hom} (W/\Lambda, \Lambda)$. 
Fixing an ``origin'' in $\Lambda(\tau)^\pitchfork$ we make 
$\Lambda_\tau(t)$ a vector function with values in ${\rm 
Hom}\,(W/\Lambda,\Lambda)$ and with the pole at $t=\tau$. 
Obviously, only free term in the expansion of this function 
to the Laurent series at $\tau$ depends on the choice of 
the ``origin'' we did to identify the affine space with the 
linear one. More precisely, the addition of a vector to the 
``origin'' results in the addition of the same vector to 
the free term in the Laurent expansion. In other words, for 
the Laurent expansion of a curve in an affine space, the 
free term of the expansion is an element of this affine 
space. Denote this element by $\Lambda^0(\tau)$. The curve 
$\tau\mapsto\Lambda^0(\tau)$ is called the {\it derivative 
curve} of $\Lambda(\cdot)$. 
 
If we restrict ourselves to the Lagrange Grassmannian 
$L(W)$, i.e. if all subspaces under consideration are 
Lagrangian w.r.t. some symplectic form $\bar\sigma$ on $W$, 
then from Remark \ref{lagr} it follows easily that the set 
$\Lambda^\pitchfork_L$ of all Lagrange subspaces 
transversal to $\Lambda$ can be considered as the affine 
space over the linear space of all self-adjoint mappings 
from $\Lambda^*$ to $\Lambda$, the velocity 
${\partial\over\partial t}\Lambda_\tau(t)$ is well defined 
self-adjoint mappings from $\Lambda^*$ to $\Lambda$, and 
the derivative curve $\Lambda^0(\cdot)$ consist of Lagrange 
subspaces. Besides if the curve $\Lambda(\cdot)$ is 
nondecreasing rank $1$ curve in $L(W)$, then 
${\partial\over\partial t}\Lambda_\tau(t)$ is a nonpositive 
definite rank $1$ self-adjoint map from $\Lambda^*$ to 
$\Lambda$ and for $t\neq \tau$ there exists a unique, up to 
the sign, vector $w(t,\tau) \in \Lambda(\tau)$ 
  such that for any $v\in \Lambda(\tau)^*$
\begin{equation}
\label{defw1} \langle v,{\partial\over\partial 
t}\Lambda_\tau(t) v\rangle= - \langle v,w(t,\tau)\rangle^2.
\end {equation} The properties of vector function $t\mapsto 
w(t,\tau)$ for a rank $1$ curve of constant weight in 
$L(W)$ can be summarized as follows ( see \cite{jac1}, 
section 7, Proposition 4 and Corollary 2): 

\begin{prop}
\label{constcor} If $\Lambda(\cdot)$ is a rank 1 curve of 
constant weight in L(W), then for any $\tau\in I$ the 
function $t\mapsto w(t,\tau)$ has a pole of order $m$ at 
$t=\tau$. Moreover, if we write down the expansion of 
$t\mapsto w(t,\tau)$ in Laurent series at $t=\tau$, 
$$w(t,\tau)=\sum_{i=1}^{m}e_i(\tau)(t-\tau)^{i-1-l}+O(1),$$ 
then the vector coefficients $e_1(\tau),\ldots,e_m(\tau)$ 
constitute a basis of the subspace $\Lambda(t)$. 
\end{prop}

The basis of vectors $e_1(\tau),\ldots, e_m(\tau)$ from the 
previous proposition is called a {\it canonical basis} of 
$\Lambda(\tau)$. Further for given $\tau$ take the 
derivative subspace $\Lambda^0(\tau)$ and let 
$f_1(\tau),\ldots, f_m(\tau)$ be a basis of 
$\Lambda^0(\tau)$ dual to the canonical basis of 
$\Lambda(\tau)$, i.e. $\bar\sigma (f_i(\tau), 
e_j(\tau))=\delta_{i,j}$. The basis 
$$(e_1(\tau),\ldots,e_m(\tau),f_1(\tau),\ldots,f_m(\tau))$$ 
of whole symplectic space $W$ is called {\it the canonical 
moving frame} of the curve $\Lambda(\cdot)$. Calculation of 
structural equation for the canonical moving frame is 
another way to obtain symplectic invariants of the curve 
$\Lambda(\cdot)$. 

For the reduced Jacobi curves of abnormal extremals of 
$(2,5)$-distribution $m$ is equal to $2$. So we restrict 
ourselves to this case. For $m=2$ the structural equation 
for the canonical moving frame has the following form (for 
the proof see \cite {jac2} Section 2, Proposition 7): 

\begin{equation}
\label{structeq} \left( 
\begin{array}{l} e_1'(\tau)\\  e_2'(\tau)\\  f_1'(\tau)\\ 
f_2'(\tau)\end{array}\right) = \left(\begin{array}{cccc} 
0&3 &0&0\\ \frac{1}{4} \rho(\tau)&0&0 &4\\ -\left(\frac 
{35}{36} A(\tau) -\frac{1}{8} 
\rho(\tau)^2+\frac{1}{16}\rho''(\tau)\right)& 
-\frac{7}{16}\rho'(\tau) &0&-\frac {1}{4}\rho(\tau) \\ 
-\frac{7}{16} \rho'(\tau) &-\frac{9}{4} \rho(\tau) & - 3 
&0\end{array}\right)\left(\begin{array}{l} 
e_1(\tau)\\e_2(\tau)\\f_1(\tau)\\f_2(\tau)\end{array}\right), 
\end{equation}
where $\rho(\tau)$ and $A(\tau)$ are the Ricci curvature 
and the density of fundamental form of the parametrized 
curve $\Lambda(\tau)$ respectively. 

Here we use the method of computation, which is slightly 
different from the method of section 3 of \cite {zelduke}, 
using again the structural equation (\ref{structeq}). 
Vector $e_1$ (and therefore $e_2$) can be found relatively 
easy (see Proposition 3.3 in \cite{zelduke}). On the other 
hand, it is difficult to compute the vectors $f_1$ and 
$f_2$, because they are defined with the help of the 
derivative curve. Instead of this one can complete the 
canonical basis $\bigl(e_1,e_2\bigr)$ of the curve somehow 
to the symplectic moving frame, find the structural 
equation for this frame and express the invariants $\rho$ 
and $A$ in terms of this structural equation. 

Indeed,
let $\tilde f_1(t),\tilde f_2(t)$ be another two vectors 
which complete the basis $\bigl(e_1(t),e_2(t)\bigr)$ to a 
symplectic moving frame in $W$. Then it is easy to show 
(see \cite{jac1}, section 7) that the structural equation 
for the frame $\bigl(e_1(t),e_2(t),\tilde f_1(t),\tilde 
f_2(t)\bigr)$ has always the form 
\begin{equation}
\label{structaux1} \left( 
\begin{array}{l} e_1'(t)\\  e_2'(t)\\ \tilde f_1'(t)\\
\tilde f_2'(t)\end{array}\right) = 
\left(\begin{array}{cccc} 0 & 3 &0&0\\ a_{21}(t)&a_{22}(t)&0& 
4\\ a_{31}(t)&a_{41}(t)&0&-a_{21}(t) \\ 
a_{41}(t)&a_{42}(t)&-3&-\alpha_{22}(t)\\ 
\end{array}\right)\left(\begin{array}{l}
e_1(t)\\ e_2(t)\\ \tilde f_1(t)\\ \tilde 
f_2(\tau)\end{array}\right). 
\end{equation}
The following lemma is the base for our method:
\begin{lemma}
\label{canlem2} The Ricci curvature $\rho(t)$ and the  
density $A(t)$ of the fundamental form 
 of the curve 
$\Lambda(t)$ satisfy

\begin{equation}
\label{auxrho} \rho=-\frac{4}{15} \Bigl(3 a_{21}+2 
a_{42}+\frac{1}{2}\alpha_{22}'+ \frac{1}{2}\alpha_{22}^2), 
\end{equation}
\begin{equation}
\label{auxA} 
A=\frac{36}{35}\Bigl(-a_{31}+\frac{9}{64}\rho^2+
\frac{1}{16}\rho''-\frac{1}{4}(a_{21})^2+\frac{1}{3} 
a'_{41}+\frac{1}{12}a''_{21}+\frac{1}{12}(a_{21} 
a_{22})'\Bigr). 
\end{equation}
\end{lemma}

{\bf Proof.} Obviously, the vectors $\tilde f_1(t), \tilde 
f_2(t)$ can be expressed by the canonical moving frame 
of $\Lambda(t)$ as follows  
\begin{equation}
\label{mu1} 
\begin{array}{c}
\tilde f_1= f_1 +\mu_{11} e_1 + \mu_{12} e_2 \\ \tilde f_2= 
f_2 +\mu_{12} e_1 + \mu_{22} e_2 
\end{array}
\end{equation}
It easily implies that 

\begin{equation}
\label{dote2} 
 \dot e_2=(a_{21}+4\mu_{12})e_1+(a_{22}+4\mu_{22})e_2+4 f_2
\end{equation}
Comparing this with the second row in 
(\ref{structeq}) we have 
\begin{eqnarray}
&~&\mu_{22}=-\frac{1}{4}a_{22} \label{mu22}\\ &~& 
a_{21}+4\mu_{12}=\frac{1}{4}\rho\label{a21} 
\end{eqnarray}
Further, using (\ref{mu22}), it is easy to obtain that 
\begin{equation}
\label{dotf2} \dot 
f_2=(a_{41}+\frac{1}{4}a_{21}a_{22}-3\mu_{11}-\mu_{12}')e_1+
(a_{42} +\frac{1}{4}(a_{22}^2+a_{22}')-6\mu_{12})e_2-3f_1 
\end{equation}
Comparing this with the forth row in 
(\ref{structeq}) we have 
\begin{eqnarray}
&~&a_{42}+\frac{1}{4}(a_{22}^2+a_{22}')-6\mu_{12}=
-\frac{9}{4}\rho\label{a42}\\ &~& 
a_{41}+\frac{1}{4}a_{21}a_{22}-3\mu_{11}-\mu_{12}'=
-\frac{7}{16}\rho'\label{a41} 
\end{eqnarray}
From equations (\ref{a21}) and (\ref{a42}) we obtain 
(\ref{auxrho}). 

Further, 
from  (\ref{a21}) and (\ref{a41}) it follows easily that 
\begin{equation}
\label{mu11} \mu_{11}=\frac{1}{3} 
a_{41}+\frac{1}{12}(a_{21}'+a_{21} a_{22})+ 
\frac{1}{8}\rho' 
\end{equation}
 Besides, by direct calculations one can get 
without difficulties that
\begin{equation}
\label{a31} \dot f_1=\Bigl(a_{31}-(\frac{1}{4}\rho+a_{21})
\mu_{12}-\mu_{11}'\Bigr)e_1-\frac{7}{16}\rho'e_2- 
\frac{\rho}{4}f_2 
\end{equation}
Comparing this with the third equation in (\ref{structeq}) 
and using (\ref{a21}) we get 
\begin{equation}
\label{preA} -\frac {35}{36} A +\frac{1}{8} 
\rho^2-\frac{1}{16}\rho''= 
 a_{31}-\frac{1}{64}\rho^2 +\frac{1}{4} a_{21}^2
 -\mu_{11}'
\end{equation}
Substituting (\ref{mu11}) in (\ref{preA}) one can easily 
obtain (\ref{auxA}). 
$\Box$ 
\medskip
 
Let us apply the previous lemma to the calculation of the 
fundamental form of (2,5)-distribution $D$. Let again 
$\{X_i\}_{i=1}^5$ be an adapted frame to $D$ and $\vec 
h=\vec h_{X_1,X_2}$ be as in (\ref{ham25}). Then it is not 
difficult to show that 
\begin{equation}
\label{vechc25} 
\begin{array}{c}
 \vec h=u_4\vec u_2-u_5\vec u_1= u_4 X_2-u_5 X_1+
\Bigl(c_{42}^4 u_4^2+(c_{42}^5-c_{41}^4)u_4 u_5-c_{41}^5 
u_5^2\Bigr)\partial_{u_4}+\\ +\Bigl(c_{52}^4 
u_4^2+(c_{52}^5-c_{51}^4)u_4 u_5-c_{51}^5 
u_5^2\Bigr)\partial_{u_5},
\end{array}
\end{equation}
where $c_{ji}^k$ are the structural functions of the frame 
$\{X_i\}_{i=1}^5$, i.e., the functions, satisfying 
$[X_i,X_j]=\sum_{k=1}^5 c_{ji}^k X_k$. 

 For any $\lambda\in (D^2)^\perp\backslash 
(D^3)^\perp$ consider the characteristic curve $\gamma$ of 
$D$ passing through $\lambda$. Under identification 
(\ref{tanid}) the reduced Jacobi curve $\widetilde 
J_\gamma$ lives in Lagrange Grassmannian $L(W_\lambda)$ of 
symplectic space $W_\lambda$, defined by (\ref{wl}). Let 
$\epsilon_1(\lambda)$ be the first vector in the canonical 
basis of the curve $t\mapsto J_\gamma(e^{t \vec h}\lambda)$ 
at the point $t=0$. 
Note that it is more convenient to work directly with 
vector fields of $(D^2)^\perp$, keeping in mind that the 
symplectic space $W_\lambda$ 
belongs to the factor space 
$T_\lambda\bigl((D^2)^\perp\bigr)/ {\rm span} \bigl(\vec 
h(\lambda),\vec e(\lambda)\bigr)$. So, in the sequel by 
$\epsilon_1(\lambda)$ we will mean both the element of 
$W_\lambda$ 
and some representative of this element in 
$T_\lambda\bigl((D^2)^\perp\bigr)$, depending smoothly on 
$\lambda$. In the last case all equalities, containing 
$\epsilon_1(\lambda)$, 
 will be assumed
modulo ${\rm span} \bigl(\vec h(\lambda),\vec 
e(\lambda)\bigr)$. By Proposition 3.4 of \cite{zelduke} the 
vector $\epsilon_1(\lambda)$ can be chosen in the form 
\begin{equation}
\label{e125eq} \epsilon_1(\lambda)= 6\Bigl( 
\gamma_4(\lambda) \partial_{u_4}+\gamma_5(\lambda) 
\partial_{u_5}\Bigr), 
\end{equation}
where 
\begin{equation}
\label{gammarel} \gamma_4 (\lambda) u_5-\gamma_5(\lambda) 
u_4\equiv 1. 
\end{equation}

Denote by $\epsilon _2(\lambda)$ the second vector in the 
canonical basis of the curve $t\mapsto J_\gamma(e^{t \vec 
h}\lambda)$. Let us compute the vector 
$\epsilon_2(\lambda)$. First note that by definition of 
vector fields $\epsilon_1$ ,$\vec e$ and relation 
(\ref{gammarel}) one has 
\begin{equation}
\label{partue} 
\partial_{u_4}=\frac{u_5}{6}\epsilon_1 ({\rm mod}\,\, \vec 
e), \,\,\partial_{u_5}=-\frac{u_4}{6}\epsilon_1 ({\rm 
mod}\,\, \vec e) 
\end{equation}
This together with (\ref{vechc25}) yields easily that 
\begin{eqnarray} 
&~& [\vec h,\partial_{u_i}]=\frac{1}{6} \alpha_i\epsilon_1 
\quad {\rm mod}\Bigl({\rm span}(\vec e,\partial 
_{u_1},\partial_{u_2},\partial_{u_3})\Bigr), \quad 
i=1,2,3,\label{hpart}\\ &~& [\vec h,\partial_{u_4}]=-\vec 
u_2 + \frac{1}{6} \alpha_4\epsilon_1 \quad {\rm 
mod}\Bigl({\rm span}(\vec e,\partial 
_{u_1},\partial_{u_2},\partial_{u_3})\Bigr), 
\label{hpart4}\\&~& [\vec h,\partial_{u_5}]=\vec u_1 + 
\frac{1}{6} \alpha_5\epsilon_1 \quad {\rm mod}\Bigl({\rm 
span}(\vec e,\partial 
_{u_1},\partial_{u_2},\partial_{u_3})\Bigr), \label{hpart5} 
\end{eqnarray} where 
\begin{equation} \label{aimal} 
\alpha_i=c_{52}^iu_4^2-(c_{42}^i+c_{51}^i)u_4u_5+ c_{41}^i 
u_5^2 
\end{equation} From the structural equation 
(\ref{structeq}), relation (\ref{e125eq}),
(\ref{gammarel}), (\ref{hpart4}),(\ref{hpart5}) one can 
easily obtain that  
\begin{equation}
\label{e225eq} \epsilon_2=\frac{1}{3}[\vec h,\epsilon_1]=
2[\vec h,\gamma_4\partial_{u_4}+\gamma_5\partial_{u_5}]=2\,X-b\,\epsilon_1
\end{equation} 
where
\begin{eqnarray}
&~&\label{gammaX}
X=\gamma_5 \vec u_1-\gamma_4 \vec u_2 +\partial_{u_3}\\
&~& b=-\frac{1}{3}\Bigl(\vec h(\gamma_4)u_5-\vec h(\gamma_5)u_4+
\alpha_4\gamma_4+\alpha_5\gamma_5\Bigr)\label{gammab}
\end{eqnarray}

 
In particular from (\ref{e225eq}),(\ref{gammaX}), and 
(\ref{gammarel}) it follows that 

\begin{eqnarray} &~&\vec 
u_1=-\frac{1}{2}(b\epsilon_1+\epsilon_2) u_4\,\, {\rm 
mod}\left({\rm span} 
 (\vec h,\partial_{u_3})\right)\nonumber\\
 &~& ~ \label{vecu1u2}\\
&~&\vec u_2=-\frac{1}{2}(b\epsilon_1+\epsilon_2)u_5\,\, 
{\rm mod}\left({\rm span} 
 (\vec h,\partial_{u_3})\right)\nonumber
\end{eqnarray} 

Let us analyze more carefully the expression for $b$.
From (\ref{gammarel}) it follows that
$$\vec h(\gamma_4)u_5-\vec h(\gamma_5)u_4=-\gamma_4\vec h(u_5)+
\gamma_5 \vec h(u_4)$$
Therefore
\begin{equation}
\label{bprom}
b=-\frac{1}{3}\left(\gamma_4\bigl(\alpha_4-\vec h(u_5)\bigr)+
\gamma_5\bigl(\alpha_5+
\vec h(u_4)\bigl)\right)
\end{equation}
On the other hand, using (\ref{vechc25}) and (\ref{aimal}), one 
has

\begin{eqnarray} &~& \alpha_5+\vec 
h(u_4)=\Bigl((c_{42}^4+c_{52}^5)u_4-(c_{41}^4+c_{51}^5)u_5 
\Bigr) u_4\nonumber 
\\ ~&~&\label{ahb}\\ 
&~& \alpha_4-\vec 
h(u_5)=-\Bigl((c_{42}^4+c_{52}^5)u_4-(c_{41}^4+c_{51}^5)u_5 
\Bigr) u_5\nonumber 
\end{eqnarray}
Substituting the last formulas into (\ref{bprom}) and using again 
(\ref{gammarel}) one has

\begin{equation}
\label{bmal}
 b= 
\frac{1}{3}\Bigl((c_{42}^4+c_{52}^5)u_4-(c_{41}^4+c_{51}^5)u_5
\Bigr) 
\end{equation}
So, the function $b$ actually does not depend on the choice 
of $\gamma_4$ and $\gamma_5$ and it is linear in $u_4$ and 
$u_5$. 

Further, one can complete the canonical basis 
$(\epsilon_1(\lambda), \epsilon_2(\lambda))$ of the 
subspace $J_\gamma(\lambda)$ to the symplectic basis in the 
symplectic space $W_\lambda$ by adding two vectors 
$\widetilde{\Phi}_1(\lambda)$ and $\widetilde 
\Phi_2(\lambda)$ such that 
\begin{eqnarray}
&~& \widetilde{\Phi}_2
=\frac{1}{2}(\vec 
u_3+u_4\partial_{u_1}+ u_5\partial_{u_2}),\label{compf2} \\
&~&~\nonumber\\
&~& \widetilde{\Phi}_1=\frac{1}{6} Y_4+b\widetilde \Phi_2,
\label{compf1} 
\end{eqnarray}
where
\begin{equation}
\label{vecYk}
 Y_4= \frac{1}{6}\Bigl( u_5\vec u_4-u_4\vec u_5+ 
\sum_{i=1}^3 \bigl(u_{5}\{u_i,u_4\}- 
u_4\{u_i,u_{5}\}\bigr)\partial_{u_i}\Bigr).
\end{equation}
Then the tuple 
\begin{equation}
\label{tupframe} \bigl(e_1(t),e_2(t),\tilde f_1(t),\tilde 
f_2(t)\bigr)\stackrel{def}{=}(e^{-t\vec 
h})_*\bigl(\epsilon_1,\epsilon_2,\widetilde {\Phi}_1, 
\widetilde{\Phi}_2\bigr)(e^{t\vec h}\lambda) 
\end{equation}
is a symplectic frame in $W_\lambda$ such that 
$(e_1(t),e_2(t))$ is canonical basis of the reduced Jacobi 
curve $J_\gamma(e^{t\vec h}\lambda)$ (in the righthand side 
of (\ref{tupframe}) one applies $(e^{-t\vec h})_*$ to each 
vector field in the brackets at the indicated point). In 
the sequel we will denote by $a_{ij}(\lambda)$ the entries 
of the matrix in the structural equation for the frame 
(\ref{tupframe}) at $t=0$ by analogy with the equation 
(\ref{structaux1}). Namely, let 

\begin{equation}
\label{hphia}
\begin{array}{l}
 \dot{\tilde f}_1(0)=[\vec h,\widetilde \Phi_1](\lambda)= 
a_{31}\epsilon_1 +a_{41} \epsilon_2-a_{21}\widetilde\Phi_2 
\\ \dot{\tilde f}_2(0)=[\vec h,\widetilde \Phi_2](\lambda)= 
a_{41}\epsilon_1 +a_{42} \epsilon_2- 3 \widetilde 
\Phi_1-a_{22}\widetilde\Phi_2 
\end{array}
\end{equation}
Let us calculate coefficient $a_{ij}$,appearing in 
(\ref{hphia}), in order to apply Lemma \ref{canlem2}.

\begin{lemma}
\label{hf2l} The following relations hold 

\begin{equation}
\label{a22f} a_{22}=-b_1-3b, 
\end{equation}
\begin{equation}
\label{a42f} a_{42}= 
-\frac{1}{4}\Pi 
\end{equation}
\begin{equation}
\label{a41f} a_{41}= 
-\frac{1}{4}\Pi b+\frac{1}{12}(\alpha_1u_4+\alpha_2u_5). 
\end{equation}

where 
 $$\Pi=(c_{32}^2 u_4-c_{31}^2 u_5)u_5-(c_{32}^1 
u_4-c_{31}^1 u_5)u_4-\Bigl(u_5\{u_3,u_4\}-u_4\{u_3,u_5\} 
\Bigr)= $$ 
\begin{equation}
\label{Dbig} 
(c_{32}^1+c_{53}^4)u_4^2+(c_{32}^2-c_{31}^1-c_{43}^4+c_{53}^5)u_4u_5- 
(c_{31}^2+c_{43}^5) u_5^2, 
\end{equation}
\begin{equation}
\label{b1def} b_1=c_{32}^3u_4-c_{31}^3u_5. \end{equation} 
\end{lemma}

{\bf Proof.} By (\ref{compf2}) and (\ref{aimal}) $$[\vec 
h,\widetilde{\Phi}_2]=\frac{1}{2}[\vec h, \vec 
u_3+u_4\partial_{u_1} +u_5\partial_{u_2}]=\frac{1}{2}([\vec 
h,\vec u_3]+ u_4[\vec h,\partial_{u_1}]+u_5[\vec 
h,\partial_{u_2}])\,\, {\rm mod} \Bigl({\rm 
span}(\partial_{u_1},\partial_{u_2})\Bigr)=$$ 
\begin{equation}
\label{calhf21} \frac{1}{2}\Bigl([\vec h,\vec 
u_3]+\frac{1}{6}(u_4 \alpha_1 +u_5 \alpha_2)\epsilon_1 
\Bigr) \quad {\rm mod} \Bigl({\rm span}(\vec h,\vec 
e,\partial_{u_1},\partial_{u_2},\partial_{u_3}) \Bigr) 
\end{equation}
Further, from definition of adapted frame (see 
(\ref{x345})) it follows that $$[\vec h,\vec u_3]=[u_4 \vec 
u_2-u_5\vec u_1,\vec u_3]= (u_4\vec u_5-u_5\vec u_4)+$$ 
$$\sum_{i=1}^3(c_{32}^i u_4-c_{31}^i u_5)\vec u_i- 
\{u_3,u_4\}\vec u_2+\{u_3,u_5\}\vec u_1$$ Using 
(\ref{compf2}),(\ref{compf2}) and (\ref{vecu1u2}), one can 
obtain from the last equation that 
\begin{equation}
\label{hu3} [\vec h,\vec u_3]=-6 \widetilde{\Phi}_1+ 
2\bigl( 
b_1+3b\bigr)\widetilde \Phi_2- \frac{1}{2} \Pi 
(b\epsilon_1+\epsilon_2) \quad {\rm mod} \Bigl({\rm 
span}(\vec h,\vec 
e,\partial_{u_1},\partial_{u_2},\partial_{u_3}) \Bigr), 
\end{equation}
where $\Pi$ and $b_1$ are as in (\ref{Dbig}) and 
(\ref{b1def}) respectively. Finally, substituting 
(\ref{hu3}) 
to (\ref{calhf21}), we get 
\begin{equation}
\label{calhf22} [\vec h,\widetilde \Phi_2]= \left( 
-\frac{1}{4}\Pi b+ 
\frac{1}{12}(\alpha_1u_4+\alpha_2u_5)\right)\epsilon_1 
-\frac{1}{4}\Pi \epsilon_2+3\widetilde{\Phi}_1+ 
(b_1+3b)\widetilde{\Phi}_2\quad {\rm mod}\Bigl({\rm 
span}(\vec h,\vec e)\Bigr), 
\end{equation}
which concludes the proof of the lemma. $\Box$ 


\begin{lemma}
\label{hF1l} 
The following relations hold
\begin{equation}
\label{a21f} a_{21}=-\Bigl( 
b\, b_1+\vec h(b)-\frac{1}{3} a_3\Bigr), 
\end{equation}
\begin{equation}
\label{a31f} a_{31}=\frac{1}{36}(\Omega-\Theta)+ 
\frac{1}{6}(\alpha_1 u_4+\alpha_2 u_5)b 
-\frac{1}{4}\Pi b^2,
\end{equation}
where 
\begin{equation}
\label{Xieq} \Theta=X_5(\alpha_4)u_4^2+\Bigl(X_5(\alpha_5)- 
X_4(\alpha_4)\Bigr)u_4 u_5-X_4(\alpha_5) u_5^2, 
\end{equation} 
\begin{equation}
\label{Xi1eq} \Omega= \sum_{i=1}^3\Bigl(u_5\{u_i,u_4\}- 
u_4\{u_i,u_5\}\Bigr)\alpha_i 
\end{equation}
\end{lemma}

{\bf Proof.} 
 By 
definition of $\widetilde {\Phi}_1$ (see (\ref{compf1}) 
\begin{equation}
\label{calF11} [\vec h,\widetilde \Phi_1]=\frac{1}{6}[\vec 
h,Y_4]+\vec h(b)\widetilde \Phi_2 +b[\vec h,\widetilde 
{\Phi}_2]. 
\end{equation} 

First let us calculate $[\vec h,Y_4]$: 
\begin{eqnarray}
&~& [\vec h,Y_4]= \left [u_4\vec u_2-u_5\vec u_1,u_5\vec 
u_4-u_4\vec u_5+ 
\sum_{i=1}^3\Bigl(u_5\{u_i,u_4\}-u_4\{u_i,u_5\}\Bigr)\partial_{u_i}\right]= 
\nonumber\\ &~&-\Bigl(\underbrace{u_4^2[\vec u_2,\vec 
u_5]-u_4u_5([\vec u_2,u_4]+[\vec u_1,\vec u_5])+u_5^2[\vec 
u_1,\vec u_4]}_{I}\Bigr)+\label{hY4}\\ &~& \underbrace{\vec 
h(u_5)\vec u_4-\vec h(u_4)\vec 
u_5}_{II}+\underbrace{\sum_{i=1}^3\Bigl(u_5\{u_i,u_4\}- 
u_4\{u_i,u_5\}\Bigr) [\vec h,\partial_{u_i}]}_{III}\,\, 
{\rm mod}\, \vec h \nonumber 
\end{eqnarray} 
Note that the vector fields $[\vec u_i,\vec u_j]$ 
calculated at the points of $(D^2)^\perp$ satisfy 

\begin{eqnarray}
 &~&[\vec u_i,\vec u_j]=\sum_{k=1}^5\Bigl(c_{ji}^k \vec u_k-
 X_k(c_{ji}^4u_4+c_{ji}^5u_5)\partial_{u_k}\Bigr)=
 \sum_{k=1}^5c_{ji}^k \vec u_k+\frac{1}{6}
 \Bigl(X_5(c_{ji}^4) 
 u_4^2+\Bigr.\nonumber\\
 &~&\Bigl.\bigl(X_5(c_{ji}^5)-X_4(c_{ji}^4)\bigr)u_4u_5-X_4(c_{ji}^5)u_5^2
 \Bigr)\epsilon_1\,\,
 {\rm mod}\Bigl({\rm span}(\partial_{u_1},
 \partial_{u_2},\partial_{u_3}),\vec e\Bigr) 
 \label{uiuj}
 \end{eqnarray}
 (the last equality here was obtained with the help of 
 (\ref{partue})).
 Substituting the last formula in the term  $I$ of
 (\ref{hY4}) one can get by direct calculation that
 \begin{equation}
 \label{IIeq}
 I=\sum_{k=1}^5\alpha_k \vec u_k+\frac{1}{6}\Theta\epsilon_1
 \quad {\rm mod}\Bigl({\rm span}(\partial_{u_1},
 \partial_{u_2},\partial_{u_3}),\vec e\Bigr), 
 \end{equation}
where $\alpha_k$ are exactly as in (\ref{aimal}) and 
$\Theta$ is defined by (\ref{Xieq}). 

Using (\ref{compf2}) and (\ref{vecu1u2}), one can obtain 
from (\ref{IIeq}) that 
\begin{eqnarray}
&~& -I+II=(\vec h(u_5)-\alpha_4)\vec u_4-(\vec 
h(u_4)+\alpha_5)\vec u_5- 2\alpha_3\widetilde 
\Phi_2+\nonumber\\ &~&~\label{I&IIeq}\\ 
&~&\frac{1}{2}(\alpha_1 u_4+\alpha_2 
u_5)(b\epsilon_1+\epsilon_2)- \frac{1}{6}\Theta 
\epsilon_1\quad{\rm mod}\Bigl({\rm span}(\partial_{u_1}, 
 \partial_{u_2},\partial_{u_3},\vec h,\vec e)\Bigr)\nonumber
\end{eqnarray}

Then from relations (\ref{ahb}),(\ref{bmal}), and 
(\ref{compf1}) it follows easily that 
\begin{eqnarray}
&~& -I+II=\left(-\frac{1}{6}\Theta+\frac{1}{2}(\alpha_1 
u_4+\alpha_2 u_5)b\right) \epsilon_1+\frac{1}{2}(\alpha_1 
u_4+\alpha_2 u_5)\epsilon_2+ \nonumber\\ 
&~&~\label{I&IIeqf}\\ &~&18 b\widetilde{\Phi}_1-(2a_3+18 
b^2)\widetilde{\Phi}_2 \quad{\rm mod}\Bigl({\rm 
span}(\partial_{u_1}, 
 \partial_{u_2},\partial_{u_3},\vec h,\vec e)\Bigr)\nonumber
\end{eqnarray} 
Note also that by (\ref{hpart}) and (\ref{Xi1eq}) the term 
$III$ of (\ref{hY4}) satisfies 
\begin{equation}
\label{IIIeq} III=\frac{1}{6}\Omega \epsilon_1 \quad{\rm 
mod}\Bigl({\rm span}(\partial_{u_1}, 
 \partial_{u_2},\partial_{u_3},\vec h,\vec e)\Bigr)
\end{equation}
Substituting (\ref{I&IIeqf}) and (\ref{IIIeq}) in 
(\ref{hY4}) and then in (\ref{calF11}) and using 
(\ref{calhf22}) one can obtain without difficulties that
\begin{eqnarray}
&~&[\vec h,\widetilde 
\Phi_1]=\left(\frac{1}{36}(\Omega-\Theta)+ 
\frac{1}{6}(\alpha_1 u_4+\alpha_2 u_5)b 
-\frac{1}{4}\Pi b^2\right)\epsilon_1+\nonumber\\ 
&~&~\label{hF1}\\ &~& \left( 
-\frac{1}{4}\Pi b+ 
\frac{1}{12}(\alpha_1u_4+\alpha_2u_5)\right)\epsilon_2+ 
\Bigl( 
b\,b_1+\vec h(b)-\frac{1}{3} a_3\Bigr)\widetilde \Phi_2, 
\nonumber 
\end{eqnarray}
which concludes the proof of the lemma. $\Box$

Substituting relations (\ref{a22f}), (\ref{a42f}), and 
(\ref{a21f}) in (\ref{auxrho}) it is not difficult to 
obtain that the Ricci curvature $\rho(\lambda)$ of the 
reduced Jacobi curve $J_\gamma(e^{t\vec h}\lambda)$ at 
$t=0$ satisfies 

\begin{equation}
\label{rhocar} 
\rho=-\frac{4}{15}(a_3-\frac{1}{2}\Pi-\frac{1}{2}\vec 
h(b_1)-\frac{9}{2}\vec 
h(b)+\frac{1}{2}b_1^2+\frac{9}{2}b^2).
\end{equation}
 
Finally substituting the last formula and formulas 
(\ref{a22f}), (\ref{a41f}), (\ref{a21f}), and (\ref{a31f}) 
in (\ref{auxA}) we obtain the following expression for the 
density $A$ of the fundamental form: \begin{eqnarray} 
\frac{35}{36} A&=& \frac{1}{36}(\Theta-\Omega)- 
\frac{1}{6}(\alpha_1 u_4+\alpha_2 u_5)b 
+\frac{1}{4}\Pi b^2+\nonumber\\ ~&~&\frac{1}{100} 
\Bigl(a_3-\frac{1}{2}\Pi-\frac{1}{2}\vec 
h(b_1)-\frac{9}{2}\vec 
h(b)+\frac{1}{2}b_1^2+\frac{9}{2}b^2\Bigr)^2-\nonumber\\ 
~&~& \frac{1}{60}\vec h\circ\vec h\Bigl( 
a_3-\frac{1}{2}\Pi-\frac{1}{2}\vec h(b_1)-\frac{9}{2}\vec 
h(b)+\frac{1}{2}b_1^2+\frac{9}{2}b^2 
\Bigr)-\label{FORMULA}\\~&~& \frac{1}{4}\Bigl( 
b\, b_1+\vec h(b)-\frac{1}{3} a_3\Bigr)^2+\frac{1}{36}\vec 
h(\alpha_1 u_4+\alpha_2 u_5)-\frac{1}{12}\vec h(\Pi 
b)-\nonumber\\~ &~& \frac{1}{12}\vec h\circ\vec h\Bigl( 
b\, b_1+\vec h(b)-\frac{1}{3} a_3\Bigr)+\frac{1}{12}\vec h
\Bigl(\bigl( 
b\, b_1+\vec h(b)-\frac{1}{3} 
a_3\bigr)(b_1+3b)\Bigr).\nonumber 
\end{eqnarray} 
 
\begin{remark}
\label{structrem} {\rm Analyzing all functions involved in 
(\ref{FORMULA}), it is not hard to see that the 
coefficients of fundamental form (which is polynomial of 
degree $4$ in $u_4$ and $u_5$) can be expressed by the 
following structural functions of the chosen adapted frame: 
\begin{eqnarray}&~& c_{3j_1}^{k_1},\,\,\,\,j_1=1,2;\,\, k_1=1,2,3; \nonumber\\
&~& c_{i_2j_2}^{k_2},\,\,\, i_2=4,5;\,\, j_2=1,2;\,\, 
k_2=1\ldots 5;\label{vazhstr}\\ &~& c_{i_3 
3}^{k_3},\,\,\,\, i_3=4,5;\,\, 
k_3=4,5.\nonumber\end{eqnarray} Note also that the 
structural functions of the third raw of (\ref{vazhstr}) 
can be expressed by structural functions of the second 
row.} $\Box$ 
\end{remark}

\begin{remark}
\label{strongb} {\rm If the chosen frame is strongly 
adapted to the distribution, then $c_{3j}^{k}=0$ for 
$j=1,2$ and therefore $$\begin{array}{c}b_1=0\\ \Pi= 
-u_5\{u_3,u_4\}+u_4\{u_3,u_5\}. 
\end{array}$$}$\Box$
\end{remark}


\section{Proof of Theorem \ref {comptheor}} 
\indent\setcounter{equation}{0}
\label{proof_sec} 

Fix some coframe $\{\om_i\}_{i=1}^5$ 
satisfying (\ref{costruct}). More precisely, take some 
section of the bundle ${\cal K}_1(M)$ over $M$. We will 
think of all $1$-forms under consideration and of the 
equations (\ref{costruct}),(\ref{bom1})-(\ref{bom3}) as the 
objects restricted on this section. Let $\{\widetilde 
X_k\}_{k=1}^5$, be the frame on $M$ dual to the coframe 
$\{\om_i\}_{i=1}^5$, namely, 
\begin{equation}
\label{dual} \om_i(\widetilde X_k)=\delta_{i,k}. 
\end{equation} 
Let \begin{equation} \label{XY} X_k=\widetilde 
X_{5-k+1},\quad 1\leq k\leq 5 \end{equation} Any frame 
$\{X_i\}_{i=1}^5$ obtained in such way will be called {\it 
Cartan's frame} of the distribution $D$. 
Note that 
\begin{equation}
\label{dYX} D(q)=span(\widetilde X _4(q),\widetilde X_5(q))=span(X_1(q), 
X_2(q)),\quad q\in M.
\end{equation}
Using (\ref{dual}) and the well-known formula 
\begin{equation}
\label{dbrack}
 d \om(V_1,V_2)=V_1 \om(V_2)-V_2\om(V_1)-\om([V_1,V_2]) 
\end{equation} 
one can find from (\ref{costruct}) the commutative 
relations for the frame $\{Y_k\}_{k=1}^5$ and therefore for 
the frame $\{X_k\}_{k=1}^5$. In particular, 
\begin{eqnarray} 
[X_1,X_2]&=&X_3+\Bigl(\bom_3(X_1)-\bom_4(X_2)\Bigr)X_1+ 
\Bigl (\bom_1(X_1)-\bom(X_2)\Bigr )X_2 \label{x12} \\ 
~[X_1,X_3] &=& 
X_4-\left(\frac{4}{3}\bom_5(X_1)+\bom_4(X_3)\right)X_1+\left(\frac{4}{3} 
\bom_6(X_1)-\bom_2(X_3)\right)X_2+\nonumber\\ ~&~& \Bigl 
(\bom_1+\bom_4\Bigr )(X_1)X_3 \label{x13} \\ 
 ~[X_2,X_3] &=& 
X_5-\left(\frac{4}{3}\bom_5(X_2)+\bom_3(X_3)\right)X_1+\left(\frac{4}{3} 
\bom_6(X_2)-\bom_1(X_3)\right)X_2+\nonumber\\ ~&~& \Bigl 
(\bom_1+\bom_4\Bigr )(X_2)X_3 \label{x23} 
 \end{eqnarray} 

These formulas together with (\ref{dYX}) imply that any 
Cartan's frame $\{X_i\}_{i=1}^5$ is adapted to the 
distribution $D$. So, starting with some Cartan's frame 
$\{X_i\}_{i=1}^5$, 
 one can  apply the formula (\ref{FORMULA}) for computation of the 
fundamental form. From (\ref{costruct}) the structural 
functions of the frame $\{X_i\}_{i=1}^5$ can be expressed 
in terms of the forms $\bom_j$ and the vector fields $X_i$. 
 
Hence our fundamental form can be expressed in terms of 
forms $\bom_j$ and vector fields $X_k$. We will compare 
this expression with expression of Cartan's tensor ${\cal 
F}$ which can be obtained by substitution of some of 
$X_k$'s into the formulas (\ref{bom1})-(\ref{bom3}). 

By Remark \ref{strongb} in order to calculate the 
fundamental form we will need 
the following commutative relations in addition to 
(\ref{x12})-(\ref{x23}):
 
\begin{eqnarray}
[X_1, X_4]&=&\Bigl (\bom_7(X_1)-\bom_4(X_4)\Bigr 
)X_1-\bom_2(X_4)X_2+ \bom_6(X_1)X_3+\nonumber\\ 
~&~&\Bigl (\bom_1+2\bom_4\Bigr 
)(X_1)X_4+\bom_2(X_1)X_5,\label{x14} 
\end{eqnarray}
\begin{eqnarray}
[X_1, X_5]&=&-\bom_4(X_5)X_1+\Bigl 
(\bom_7(X_1)-\bom_2(X_5)\Bigr )X_2+\bom_5(X_1)X_3+ 
\nonumber\\ 
~&~&\bom_3(X_1)X_4+\Bigl (2\bom_1+\bom_4\Bigr 
)(X_1)X_5,\label{x15}\end{eqnarray} 
\begin{eqnarray}
[X_2, X_4]&=&\Bigl (\bom_7(X_2)-\bom_3(X_4)\Bigr 
)X_1-\bom_1(X_4)X_2+\bom_6(X_2)X_3+ 
 \nonumber\\
~&~&\Bigl (\bom_1+2\bom_4\Bigr 
)(X_2)X_4+\bom_2(X_2)X_5,\label{x24}
\end{eqnarray}
\begin{eqnarray}
[X_2, X_5]&=&-\bom_3(X_5)X_1+\Bigl 
(\bom_7(X_2)-\bom_1(X_5)\Bigr )X_2+ \bom_5(X_2)X_3+ 
 \nonumber\\ 
~&~& \bom_3(X_2)X_4+\Bigl (2\bom_1+\bom_4\Bigr 
)(X_2)X_5,\label{x25} 
\end{eqnarray}
\begin{equation}
[X_3, X_4]=
\bigl(\bom_1+2\bom_4\bigr)(X_3)X_4+\bom_2(X_3)X_5\,\,\,{\rm 
mod}\Bigl(\rm{span}(X_1,X_2,X_3)\Bigr),\label{x34} 
\end{equation}
\begin{equation}
[X_3, X_5]=
\bom_3(X_3)X_4 + \bigl(2\bom_1+\bom_4\bigr)(X_3)X_5 
,\,\,{\rm mod}\Bigl(\rm{span}(X_1,X_2,X_3)\Bigr). 
\label{x35} 
\end{equation}
These relations can be easily obtained by applying formulas 
(\ref{dual}) and (\ref{dbrack}) to (\ref{costruct}). 

First note that for Cartan's frame from commutative 
relations (\ref{x14})-(\ref{x25}) it follows that 
\begin{equation}
\label{funbom} b=(\bom_1+\bom_4)(X_2) 
u_4-(\bom_1+\bom_4)(X_1) u_5. 
\end{equation} 
Further, from (\ref{b1def}), (\ref{x13}), (\ref{x23}), and 
(\ref{funbom}) it follows that for Cartan's frames
\begin{equation}
\label{c123b} b_1
=b 
\end{equation}
In addition, substituting (\ref{x13})-(\ref{x25}) in 
(\ref{aimal}) and (\ref{Dbig}) one can obtain by direct 
computations that for Cartan's frames 
 
\begin{equation}
\label{Da3} \Pi=-\frac{4}{3} \alpha_3 
\end{equation}

After substitution of (\ref{c123b}) and (\ref{Da3}) in 
(\ref{FORMULA}) most of the terms of (\ref{FORMULA}) are 
surprisingly canceled and we obtain that in Cartan's frame 
the density $A$ of fundamental form satisfies
\begin{equation}
\label{fundAexp} A=\frac{1}{35}\Bigl(\Theta+ \vec 
h(\alpha_1 u_4+\alpha_2 u_5)-\Omega-6 (\alpha_1 
u_4+\alpha_2 u_5)b\Bigr) 
\end{equation}

Let us try to express the density $A$ in the most 
convenient way in terms of forms $\bom_j$ and vector fields 
$X_k$. For this first note that the term 
 $\vec 
h(\alpha_1 u_4+\alpha_2 u_5)$ can be written in the 
following form: 
\begin{eqnarray}
 \vec h(\alpha_1 u_4+\alpha_2 u_5)&=&\Theta_1+\vec 
h(u_4)(\alpha_1+u_4{\partial\over\partial 
u_4}\alpha_1+u_5{\partial\over\partial u_4}\alpha_2) 
+\nonumber\\ ~&~&\vec h(u_5)(\alpha_2+u_4 
{\partial\over\partial 
u_5}\alpha_1+u_5{\partial\over\partial 
u_5}\alpha_2),\label{ha1u4a2u5} 
 \end{eqnarray} 
where \begin{equation} \label{theta1} \Theta_1=
X_2(\alpha_1)u_4^2+\Bigl(X_2(\alpha_2)-X_1(\alpha_1)\Bigr)u_4u_5-X_1(\alpha_2) 
u_5^2.\end{equation} 

Note that $\Theta$ and $\Theta_1$ are the only terms of $35 
A$ containing expression of the form 
 $X_k\bom_i(X_j)$. Let us analyze the sum 
 $\Theta+\Theta_1$.
For this let us introduce some notations. Given 1-form 
$\omega$ and two vector fields denote by $ {\cal 
S}(\om,V_1,V_2)$ the following expression 
\begin{equation} \label{Sop} 
 {\cal S}(\om,V_1,V_2)\stackrel{def}{=}V_1\om(V_2)-V_2\om(V_1)
 \end{equation}
 Denote also by ${{\cal W}_{u_4,u_5}}$ the following $1$-form on 
$M$: 
\begin{equation}
\label{thetas} 
  {{\cal W}_{u_4,u_5}}=u_4^2\bom_3+u_4u_5(\bom_1-\bom_4)-u_5^2\bom_2
 \end{equation} 
Then from (\ref{aimal}), (\ref{Xieq}),(\ref{theta1}), and 
commutative relations (\ref{x14})-(\ref{x25}) one can 
obtain by direct calculations that 
\begin{eqnarray}
 \Theta+\Theta_1&=& {\cal S}({\cal W}_{u_4,u_5},X_5,X_2) u_4^2+ \Bigl({\cal S}
 ({\cal W}_{u_4,u_5},X_1,X_5)+\nonumber\\~&~& {\cal S}({\cal W}_{u_4,u_5},X_2,X_4)\Bigr)u_4u_5+
 {\cal S}({\cal W}_{u_4,u_5},X_4,X_1)u_5^2 
\label{sumtheta}\end{eqnarray} 

The "commutative" nature of $\Theta+\Theta_1$ suggests an 
idea to compare the polynomial $35 A$ with the following 
polynomial of degree 4 in $u_4$ and $u_5$: 
\begin{eqnarray}
{\cal B}&\stackrel{def}{=}&d {\cal W}_{u_4,u_5}(X_5,X_2) 
u_4^2+ \Bigl(d{\cal W}_{u_4,u_5}(X_1,X_5)+\nonumber\\~&~& 
d{\cal W}_{u_4,u_5}(X_2,X_4)\Bigr)u_4u_5+ 
 d{\cal W}_{u_4,u_5}(X_4,X_1)u_5^2\label{Bpol}
\end{eqnarray} 
 (the polynomial ${\cal B}$ is obtained from the righthand side of 
 (\ref{sumtheta}) by replacing the operation ${\cal S}$ with the 
 operation of exterior differential $d$).
Substituting expressions (\ref{vechc25}) for $\vec h$, 
(\ref{aimal}) for $\alpha_i$, (\ref{funbom}) for $b$, 
formulas (\ref{ha1u4a2u5}) and (\ref{sumtheta}) in 
(\ref{fundAexp}), then formulas (\ref{thetas}) and identity 
(\ref{dbrack}) in (\ref{Bpol}) and using commutative 
relations (\ref{x14})-(\ref{x25}), one can get by long but 
direct calculations the following 
\begin{eqnarray}
35A-{\cal B}&=& \Xi_{u_4,u_5}(X_5,X_2) u_4^2+ \Bigl(\Xi 
_{u_4,u_5}(X_1,X_5)+\nonumber\\~&~& \Xi 
_{u_4,u_5}(X_2,X_4)\Bigr)u_4u_5+ 
 \Xi_{u_4,u_5}(X_4,X_1)u_5^2\label{diffpol},
\end{eqnarray} 
where $\Xi_{u_4,u_5}$ is the following 2-form on $M$: 
\begin{equation}
\label{2form} 
\Xi_{u_4,u_5}=u_4^2\bom_3\wedge(\bom_1-\bom_4)-2u_4u_5\bom_3\wedge\bom_2+u_5^2\bom_2
\wedge(\bom_1-\bom_4). \end{equation} 

On the other hand, from relations (\ref{bom1})-(\ref{bom3}) 
and duality (see relations (\ref{dual}) and (\ref{XY})) it 
follows without difficulties that 
\begin{eqnarray}
&~&{\cal B}+\Xi_{u_4,u_5}(X_5,X_2) u_4^2+ \Bigl(\Xi 
_{u_4,u_5}(X_1,X_5)+\Xi _{u_4,u_5}(X_2,X_4)\Bigr)u_4u_5+ 
\nonumber\\~&~& 
 \Xi_{u_4,u_5}(X_4,X_1)u_5^2=-\Bigl(A_1 u_4^4-4 A_2u_4^3u_5 +6 
A_3u_4^2u_5^2-4 A_4u_4u_5^3+A_5u_5^4 \Bigr)\label{bomcol}, 
\end{eqnarray}
where coefficients $A_i$ are exactly as in (\ref{cartf}). 
So, \begin{equation} \label{Aai} 35A= -\Bigl(A_1 u_4^4-4 
A_2u_4^3u_5 +6 A_3u_4^2u_5^2-4 A_4u_4u_5^3+A_5u_5^4 \Bigr)
\end{equation} 
Finally, we can find the tangential fundamental form 
$\AA_q$.
For given $v\in D(q)$ by duality one has $v=\om_4(v) 
X_2+\om_5(v)X_1$. Take $\lambda\in (D^2)^\perp(q)$ with 
\begin{equation}
\label{uom} u_4=\om_4(v),\quad u_5=-\om_5(v).
\end{equation}
 Let $\vec 
h=\vec h_{X_1,X_2}$ as in (\ref{ham25}). Then by 
construction $\pi_* \vec h(\lambda)=v$. Further from 
(\ref{tangfund}) and (\ref{densX12}) 
$$35\AA_q(v)=35 {\cal A}_\lambda(\vec h(\lambda))=35 
A(\lambda)$$ 
In order to calculate $A(\lambda)$ one substitutes
(\ref{uom}) in (\ref{Aai}). Comparing the result of this 
substitution with (\ref{cartf}), one obtains (\ref{FAA}), 
which concludes the proof of Theorem \ref{comptheor}.


\end{document}